\magnification=\magstep1 
\overfullrule=0pt      \input epsf 
 \def\leaderfill{\leaders\hbox to 1em{\hss.\hss}\hfill} 
 \def\la{\lambda}   
  \def\i{{\rm i}} 
\def\si{\sigma} \def\eps{\epsilon}    
 
\font\huge=cmr10 scaled \magstep2

\font\smcap=cmcsc10     \font\smit=cmmi7  \font\smal=cmr7
  \def\G{\Gamma}  
\input amssym.def
\def\Z{{\Bbb Z}} \def\R{{\Bbb R}} \def\Q{{\Bbb Q}}  
\def\C{{\Bbb C}}   \def\M{{\Bbb M}}
\font\smit=cmmi7     \def\H{{\Bbb H}}   \def\g{{\frak g}}

\centerline{{\huge Monstrous Moonshine: The first twenty-five years}}
\bigskip

\centerline{{\smcap Terry Gannon}}\medskip

\centerline{{Abstract}}

\medskip\noindent{{\smal Twenty-five years ago, Conway and Norton
published in this journal\footnote{$^\dagger$}{{\smal 
This is an invited paper for the {\smit Bull}.\ {\smit London}\
{\smit Math}.\ {\smit Soc}. Comments are welcome.}}
 their remarkable paper `Monstrous Moonshine',
proposing a completely unexpected relationship between finite simple groups
and modular functions. This paper reviews the progress made in broadening and
understanding that relationship.}}\footnote{}{{\smal 2000 {\smit Mathematics}\
{\smit Subject}\ {\smit Classification}\ 11F22, 17B69, 17B67, 81T40}}

\bigskip\centerline{\it 1. Introduction}
\medskip

\noindent It has been approximately twenty-five years since John McKay remarked that
 $$196\,884=196\,883+1\ .\eqno(1.1)$$
That time has seen the discovery of important structures, the establishment of
another deep connection between number theory and algebra, and a reinforcement
of a new era of cooperation between pure mathematics and mathematical
physics. It is a beautiful and accessible example of how mathematics can be driven by strictly conceptual
concerns, and of how the particular and the general can feed off each other. 
Now, six years after Borcherds' Fields Medal, 
the original flurry of activity is over; the new period should be one
of consolidation and generalisation and should witness the gradual movement of 
this still rather esoteric corner of mathematics toward the mainstream.  

The central question McKay's equation (1.1) raises, is: What does the $j$-function
(the left side) have to do with the Monster finite group (the right side)?
Many would argue that we still don't have our finger on the essence of the matter.
But what is clear is that we understand far more about this central question
today than we did in 1978. Today we say that there is a vertex operator
algebra, called the {\it Moonshine module} $V^\natural$, which interpolates
between the left and right sides of (1.1): its automorphism
group is the Monster and its graded dimension is the $j$-function ($-744$).

This paper tries to summarise this work of the past twenty-five
years in about as many pages. The original article [{\bf 24}] is 
still very readable and contains a wealth of information not found in other 
sources. Other reviews are [{\bf 21}], [{\bf 87}], [{\bf 12}],
[{\bf 39}], [{\bf 90}], [{\bf 48}], [{\bf 15}], [{\bf 78}], [{\bf 102}], [{\bf 16}], 
[{\bf 46}] and the introductory
chapter in [{\bf 44}], and each has its own emphasis. Our own bias here has been to
breadth at the expense of depth, which probably limits this review to be a mere annotated
sampling of representative literature.

\bigskip\centerline{{\it 2. Background}}\medskip

In \S2.1 we describe the finite simple groups and in
particular the Monster. In \S2.2 we focus on the modular groups and
functions which arise in Monstrous Moonshine.

\medskip{{\it 2.1. The Monster.}}
By definition, a simple group is one whose only normal subgroups are the
trivial ones: $\{1\}$ and the group itself. The importance of the
finite simple groups  lies in their role as building blocks,
in the sense that any finite group $G$ can be constructed from $\{1\}$ by
extending successively by a (unique up to order) sequence of finite simple
groups. For example the symmetric group $S_4$ arises in this way from the
cyclic groups $C_2,C_2,C_3,C_2$.

A formidable accomplishment last century was determining the explicit list
of all finite simple groups. See e.g.\ [{\bf 49}] for more details and references.
These groups are:\smallskip

\item{(i)} the cyclic groups $C_p$, $p$ prime;

\item{(ii)} the alternating groups $A_n$, $n\ge 5$;

\item{(iii)} 16 infinite families of groups of Lie type;

\item{(iv)} 26 sporadic groups.\smallskip

An example of a family of finite simple group of Lie type is PSL$_n({\Bbb F}_q)$, i.e.\
the group of $n\times n$ matrices with determinant 1, and entries from the
finite field ${\Bbb F}_q$, quotiented out by its centre
(the scalar matrices $aI$, where $a^n=1$).

The mysteriousness of the sporadics is due to their falling outside those infinite 
families. They range in size from the Mathieu group $M_{11}$, with order 7920 
and discovered in 1861, to the Monster ${\Bbb M}$, with order
$$|\M|=2^{45}\cdot 3^{20}\cdot 5^9\cdot 7^6\cdot 11^2\cdot 13^3\cdot 17\cdot
19\cdot  23\cdot  29\cdot  31\cdot  41 \cdot 47 \cdot 59 \cdot 71
\approx 8\times 10^{53}\ .\eqno(2.1)$$
The existence of $\M$ was conjectured in 1973 by Fischer and Griess, and
finally constructed in 1980 by Griess [{\bf 50}]. Most sporadics arise in $\M$
(e.g.\ as quotients of subgroups). We'll encounter many of these sporadics in the
coming pages, but most of our attention will be directed at $\M$.

Griess showed in fact that $\M$ was the automorphism group of a 196883-dimensional
commutative nonassociative algebra, now called the {\it Griess algebra},
but the construction was somewhat artificial. We now understand [{\bf 44}] the
Griess algebra
as the first nontrivial tier of an infinite-dimensional graded algebra, the {\it Moonshine module}
$V^\natural$, which lies at the heart of Monstrous Moonshine. We'll discuss
$V^\natural$ in \S4.2; we will find that it
has a very rich algebraic structure, is conjectured to obey a strong
uniqueness property, and has automorphism group $\M$.

The Monster $\M$ has a remarkably simple presentation. As with any noncyclic
finite simple group, it is generated by its involutions (i.e.\ elements of
order 2) and so will be a homomorphic image of a Coxeter group. Let ${\cal G}_{pqr}$,
$p\ge q\ge r\ge 2$, be the graph consisting of three strands of lengths
$p+1,q+1,r+1$, sharing a common endpoint. Label the $p+q+r+1$ nodes as in
Figure 1. 
 Given any graph ${\cal G}_{pqr}$,
define $Y_{pqr}$ to be the group consisting of a generator for each node,
obeying the usual Coxeter group relations (i.e.\ all generators are involutions,
and the product $gg'$ of two generators has order 3 or 2, depending on whether
or not the two nodes  are adjacent), together with one more relation:
$$(ab_1b_2ac_1c_2ad_1d_2)^{10}=1\ .\eqno(2.2)$$
The groups $Y_{pqr}$, for $p\le 5$, have now all been identified (see e.g.\
[{\bf 61}]).
Conway conjectured and, building on work by Ivanov [{\bf 60}], Norton proved [{\bf 98}]
that $Y_{555}\cong Y_{444}$ is the `Bimonster', the wreathed-square ${\Bbb M}
\wr C_2$ of the Monster (so has order $2|{\Bbb M}|^2$). A closely related presentation
of the Bimonster has 26 involutions as generators and has relations given by
the incidence graph of the 
projective plane of order 3; the Monster itself arises from 21 involutions and
the affine plane of order 3 [{\bf 25}]. Likewise, $Y_{553}\cong
Y_{443}\cong \M\times C_2$. Other sporadics
arise in e.g.\ $Y_{533}\cong Y_{433}$ (the Baby Monster ${\Bbb B}$), $Y_{552}\cong Y_{442}$
(the Fischer group $Fi_{24}'$), and $Y_{532}\cong Y_{432}$ (the Fischer group
$Fi_{23}$). The Coxeter groups of ${\cal G}_{555}$, ${\cal G}_{553}$, ${\cal G}_{
533}$, ${\cal G}_{ 552}$, and ${\cal G}_{ 532},$ are all infinite groups of
hyperbolic reflections in e.g.\ $\R^{17,1}$, and contain
copies of groups such as the affine $E_8$ Weyl group, so the geometry
here should be quite pretty. What role, if any,
these remarkable presentations have in Moonshine hasn't been established yet. 
As a first step though, [{\bf 97}] identifies in Aut$(V^\natural)$ the 21 
involutions generating $\M$.

\medskip\epsfysize=1.5in\centerline{ \epsffile{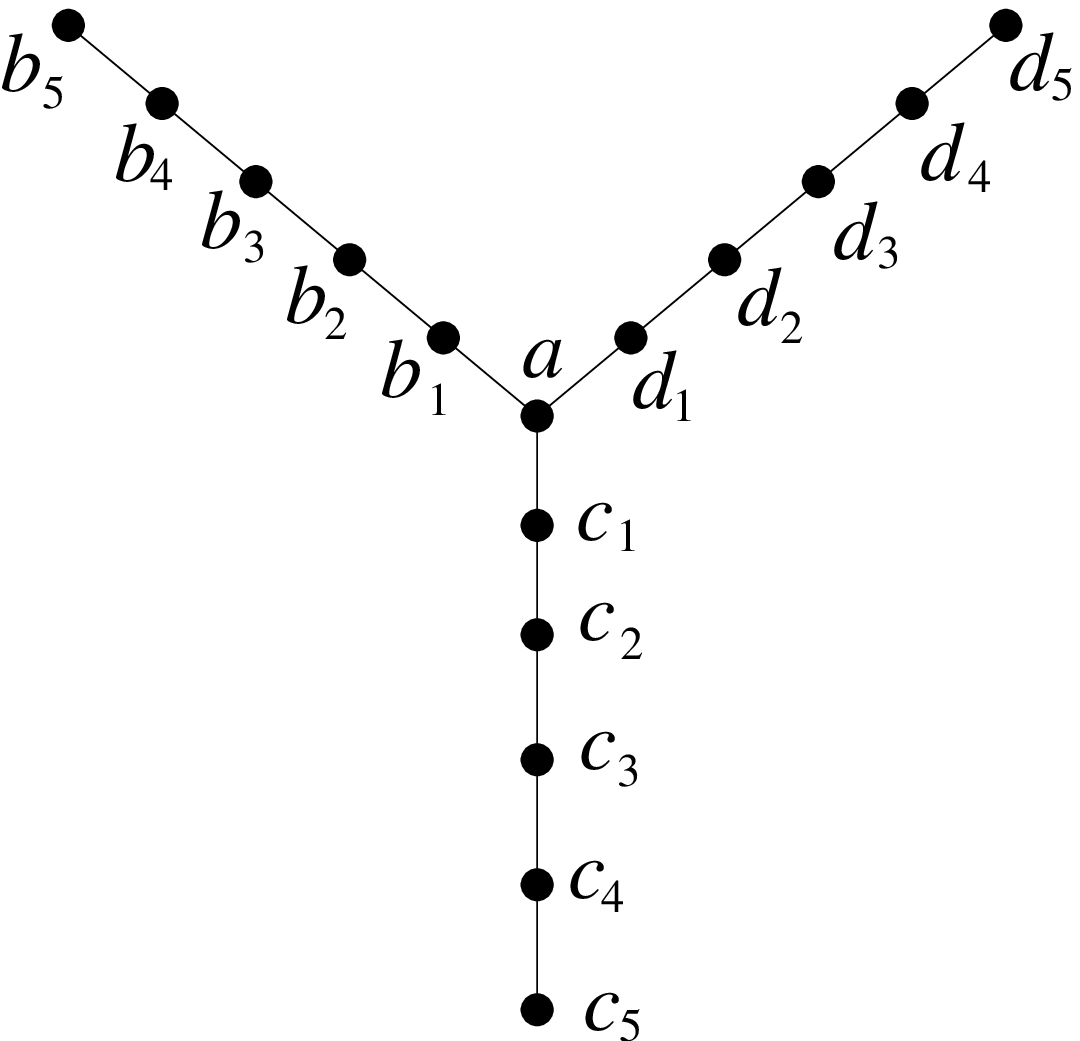}}\medskip
\centerline{{\bf Figure 1. The graph ${\cal G}_{555}$ presenting the Bimonster}}\medskip

The Monster has 194 conjugacy classes, and so that number of irreducible
representations. Its character table (and much other useful information)
is given in the Atlas [{\bf 22}], where we also find analogous data for
the other simple groups of `small' order. For example,
we find that $\M$ has exactly 2, 3, 4 conjugacy classes of elements of
order 2, 3, 4, respectively --- these classes are named 2A, 2B, 3A, etc.
We also find that the dimensions of the smallest irreducible representations
of $\M$ are 1, 196883, 21296876, and 842609326. This is the same 196883 as
on the right side of  (1.1), and as the dimension of the Griess algebra.

\medskip{{\it 2.2. The $j$-function.}}
The group SL$_2(\R)$, consisting of $2\times 2$ matrices of determinant 1 with
real entries, acts on the upper half-plane $\H:=\{\tau\in\C\,|\,{\rm Im}(\tau)>0\}$
by fractional linear transformations
$$\left(\matrix{a&b\cr c&d}\right).\tau={a\tau+b\over c\tau+d}\ .\eqno(2.3)$$
Of course this is really an action of PSL$_2(\R):={\rm SL}_2(\R)/\{\pm I\}$
on $\H$, but it is more convenient to work with SL$_2(\R)$. $\H$ is the
hyperbolic plane, one of the three possible geometries in two dimensions (the others
are the sphere and the Euclidean plane), and PSL$_2(\R)$ is its group of
orientation-preserving isometries.

Let $G$ be a discrete subgroup of SL$_2(\R)$. Then the space $G\backslash\H$
has a natural structure of an orientable surface, and inherits a complex
structure from $\H$ (so can be regarded as a complex curve). By the {\it genus}
of the group $G$, we mean the genus of the resulting real surface $G\backslash\H$.
For example, the choice $G={\rm SL}_2(\Z)$ yields 
the sphere with one puncture, so SL$_2(\Z)$ has genus 0. Moreover, any curve
$\Sigma$ with genus $g$ and $n$ punctures, for $3g+n>3$, is equivalent as a
complex curve to the space $G\backslash\H$, for some subgroup $G$ of SL$_2(\R)$
isomorphic to the fundamental group $\pi_1(\Sigma)$.

The most important choice for $G$ is SL$_2(\Z)$, thanks
to its interpretation as the modular group of the torus. Most groups $G$ of
interest are commensurable with SL$_2(\Z)$, i.e.\ $G\cap{\rm SL}_2(\Z)$
has finite index in both $G$ and SL$_2(\Z)$. Examples of these
are the congruence subgroups
$$\eqalignno{\Gamma(N):=&\,\{\left(\matrix{a&b\cr c&d}\right)\in{\rm SL}_2(\Z)\,|\,
\left(\matrix{a&b\cr c&d}\right)\equiv \left(\matrix{1&0\cr 0&1}\right)\
({\rm mod}\ N)\}\ ,&(2.4a)\cr
\Gamma_0(N):=&\,\{\left(\matrix{a&b\cr c&d}\right)\in{\rm SL}_2(\Z)\,|\,
N\ {\rm divides}\ c\}\ .&(2.4b)}$$
For example $\Gamma_0(N)$ has genus 0 for $N=2,13,25$, while $N=50$ has genus 2 and
$N=24$ has genus 3.
The following definition includes all groups arising in Monstrous Moonshine.

\medskip{{\smcap Definition 1.}} {\it Call a discrete subgroup $G$ of
SL$_2(\R)$ a} moonshine-type modular group, {\it if it contains some $\Gamma_0(N)$, and
also obeys the condition that}
$$\left(\matrix{1&t\cr 0&1}\right)\in G\ {\rm iff}\ t\in\Z\ .$$

Such a modular group is necessarily commensurable with SL$_2(\Z)$. Note that
for such a $G$, any  meromorphic function $f:G\backslash\H\rightarrow\C$ will
have a Fourier expansion of the form $f(\tau)=\sum_{n=-\infty}^\infty a_nq^n$,
where $q=e^{2\pi\i\tau}$.

\medskip{{\smcap Definition 2.}} {\it Let $G$ be any subgroup of
SL$_2(\R)$ commensurable with SL$_2(\Z)$. By a} modular function $f$ for $G$
{\it we mean any meromorphic function $f:\H\rightarrow\C$, such that}
$$f({a\tau+b\over c\tau+d})=f(\tau)\qquad \forall\left(\matrix{a&b\cr c&d}
\right)\in G\eqno(2.5)$$
{\it and such that, for any $A\in{\rm SL}_2(\Z)$, the function $f(A.\tau)$ has
Fourier expansion of the form $\sum_{n=-\infty}^\infty b_n q^{n/N}$ for some
$N$ and $b_n$ (both depending on $A$), and where $b_n=0$ for all but finitely
many negative $n$.}\medskip

This definition simply states that $f$ is a meromorphic function 
on the {\it compact} surface $\Sigma_G:=
G\backslash\overline{\H}$, where $\overline{\H}:=\H\cup\Q\cup\{\i\infty\}$.
The $G$-orbits of $\Q\cup\{\i\infty\}$ are called {\it cusps}; their role
is to fill in the punctures of $G\backslash\H$, compactifying the surface, as
there are much fewer meromorphic functions on compact surfaces than
on noncompact ones (compare the Riemann sphere to the complex plane!).

We are especially interested  in genus 0 groups $G$ of moonshine-type. Their modular
functions are particularly easy to characterise: there will be a unique modular
function $J_G$ for $G$, with $q$-expansion of the form
$$J_G(\tau)=q^{-1}+\sum_{n=1}^\infty a_n q^n\ ;\eqno(2.6)$$
the modular functions for $G$ are precisely the rational functions
$f(\tau)={{\rm poly}(J_G(\tau))\over {\rm poly}(J_G(\tau))}$ in $J_G$.
This function $J_G$ is called the (normalised) {\it Hauptmodul} for the genus 0 group $G$.
For example, the modular group SL$_2(\Z)$ has Hauptmodul
$$J_{{\rm SL}_2(\Z)}(\tau)=J(\tau)=q^{-1}+196884\,q+214\,93760\,q^2+
8642\,909970\,q^3+\cdots\ .\eqno(2.7)$$
This 196884 is the same as that on the left-side of (1.1). Historically, in place
of this Hauptmodul was the equivalent
$$j(\tau)={(\theta_2(\tau)^8+\theta_3(\tau)^8+\theta_4(\tau)^8)^3\over 8\eta(\tau)^{24}}
=J(\tau)+744\ .$$

As we know, there are other genus 0 modular groups. For example the Hauptmoduls for
$\Gamma_0(2)$, $\Gamma_0(13)$, and $\Gamma_0(25)$, are respectively
$$\eqalignno{J_2(\tau)=&\,q^{-1}+276q-2048q^2+ 11202q^3 -49152q^4+ 184024q^5+\cdots\ ,&(2.8a)\cr
J_{13}(\tau)=&\,q^{-1}-q+ 2q^2+ q^3 +2q^4 -2q^5 -2q^7 -2q^8 +q^9+\cdots\ ,&(2.8b)\cr
J_{25}(\tau)=&\,q^{-1}-q+q^4 +q^6 -q^{11} -q^{14} +q^{21}+q^{24}-q^{26}+\cdots\ .&(2.8c)}$$

Thompson [{\bf 115}] proved there are only finitely many modular groups of moonshine-type
in each genus. Cummins [{\bf 28}] has found all of these of genus 0 and 1. In particular
there are precisely 6486 genus 0 moonshine-type groups. Exactly 616 of these
have Hauptmoduls with rational (in fact integral) coefficients, the remainder have
cyclotomic integer coefficients. There are some
natural equivalences (e.g.\ a Galois action) which collapse this number to
371, 310 of which have integral Hauptmoduls.

In genus $>0$, two functions are needed to generate the function field. A
complication facing the development of a higher-genus Moonshine is that,
unlike the situation in genus 0 considered here, there is no canonical choice 
for these generators. 

See e.g.\ [{\bf 92}] for a very readable account of some of the
circle of ideas meandering through this subsection. Modular functions are 
discussed in e.g.\ [{\bf 75}].

\bigskip\centerline{{\it 3. The Monstrous Moonshine conjectures}}\medskip

 The number on
the left of (1.1) is the first nontrivial coefficient of the $j$-function, and
the numbers on the right are the dimensions of the smallest irreducible
representations of the Fischer--Griess Monster ${\Bbb M}$. On the one
side we have a modular function; on the other, a sporadic finite
simple group. Moonshine is the explanation and generalisation of this
unlikely connection.

But first, why can't (1.1) merely be a coincidence? This is soon dispelled
by comparing the next few coefficients of $J$ with the dimensions of 
irreducible representations of $\M$:
$$\eqalignno{214\,93760=&\,212\,96876+196883+1\ ,&(3.1a)\cr
8642\,99970=&\,8426\,09326+212\,96876+2\cdot 196883+2\cdot 1\ .&(3.1b)}$$

\medskip{{\it 3.1. The fundamental conjecture of Conway and Norton.}}
The central structure in the attempt to understand equations (1.1) and (3.1)
 is an infinite-dimensional graded module for the Monster:
$$V=V_0\oplus V_1\oplus V_2\oplus V_3\oplus\cdots\ .\eqno(3.2a)$$
If we let $\rho_d$ denote the $d$-dimensional
irreducible representation of $\M$, then the first few subspaces will be
$V_{0}=\rho_1$, $V_1=\{0\}$,
$V_2=\rho_1\oplus \rho_{196883}$, and 
$V_3=\rho_1\oplus\rho_{196883}\oplus\rho_{21296876}$. This module is to have
 graded dimension
$${\rm dim}_V(\tau)=\sum_{n=0}q^n{\rm dim}(V_n)=1+196884q^2+214\,93760q^3+\cdots
=q J(\tau)\ .\eqno(3.2b)$$
Of course, (3.2b) alone certainly doesn't uniquely determine $V$, but assume
for now this $V$ has been found. Thompson [{\bf 114}] suggested studying 
in addition the graded traces
$$T_g(\tau):=q^{-1}\sum_{n=0}^\infty {\rm ch}_{V_n}(g)\,q^n\eqno(3.2c)$$
for all $g\in\M$, where the ch$_{V_n}$ are characters. As taking $g=1$ recovers
$J$, (3.2c) is
a natural twist of (3.2b). The functions $T_g$ are now called the {\it McKay--Thompson
series}.

Conway and Norton conjectured [{\bf 24}] that for each element $g$ of the Monster $\M$,
$T_g$ is the Hauptmodul
$$J_{G_g}(\tau)=q^{-1}+\sum_{n=1}^\infty a_n(g)\,q^n\eqno(3.3)$$
for a genus 0 subgroup $G_g$ of SL$_2(\R)$. 
So for each $n$ the coefficient $g\mapsto a_n(g)$ defines a character ch$_{V_n}(g)$ of $\M$. 
They explicitly identify
each of the groups $G_g$; these groups each contain $\G_0(N)$ as a normal subgroup, for
some  $N$ dividing $o(g)\,{\rm
gcd}(24,o(g))$ ($o(g)$ is the order of $g$), and the quotient group $G_g/\Gamma_0(N)$
has exponent 2 (or 1). 

Since $T_g=T_{hgh^{-1}}$ by definition, there are at most 194 distinct McKay--Thompson
series. 
All coefficients $a_n(g)$ are integers (as are in fact most
entries of the character table of $\M$). 
This implies that  $T_{g}=T_{h}$ whenever
the cyclic subgroups $\langle g\rangle$ and $\langle h\rangle$ are equal. 
In fact, the total number of distinct McKay--Thompson series
$T_g$ arising in Monstrous Moonshine turns out to be only 171. 
The first 50 coefficients $a_n(g)$ of each $T_g$ are given in [{\bf 91}].
Together with the recursions given in
\S3.3 below, this allows one to effectively compute arbitrarily many coefficients $a_n(g)$
of the Hauptmoduls. It is also this which uniquely defines $V$, up to equivalence,
as a graded $\M$-module.

For example, there are two different conjugacy classes of
order 2 elements.  One of these gives the Hauptmodul $J_2$ in (2.8a), while
the other corresponds to (3.4) below.
Similarly,  (2.8b) corresponds to an order 13 element, but $J_{25}$ in (2.8c) 
doesn't equal any $T_g$. Recall that there are exactly 616 Hauptmoduls of
moonshine-type with integer coefficients, so most of these don't arise as
$T_g$. Recently [{\bf 23}], a fairly simple characterisation has been found of the 
groups  arising as $G_g$ in Monstrous Moonshine. Their proof of this
characterisation is by exhaustion.

Conway coined this conjecture {\it Monstrous Moonshine}. The word `moonshine'
here is English slang for `insubstantial or unreal', `idle talk or speculation',
`an illusive shadow'. It is meant to give the impression that matters here
are dimly lit, and that [{\bf 24}] is `distilling information illegally' from
the character table of $\M$.

Monstrous Moonshine  began, unofficially, in 1975 when Andrew Ogg remarked
that the list of primes $p$ for which the group
$$\Gamma_0(p)+:=\langle \Gamma_0(p),{1\over \sqrt{p}}\left(\matrix{0&-1\cr p&0}
\right)\rangle\eqno(3.4)$$
has genus 0, is precisely equal to the list of primes $p$ dividing the order
of $\M$. Indeed, in the tables of [{\bf 24}] we find that, for each prime $p$
dividing $|\M|$,
an element $g$ of $\M$ of order $p$ is assigned the group $G_g=\Gamma_0(p)+$.

\medskip{{\it 3.2. Lie theory and Moonshine.}}
McKay not only noticed (1.1), but also observed that
$$j(\tau)^{{1\over 3}}=q^{-{1\over 3}}\,(1+248q+4124q^2+34\,752q^3+\cdots)\ .\eqno(3.5)$$
The point is that 248 is the dimension of the defining representation of the $E_8$  Lie
group, while $4124=3875+248+1$ and $34\,752=30\,380+3875+2\cdot 248+1$.
Incidentally, $j^{{1\over 3}}$ is a generating modular function for the genus-0
group $\Gamma(3)$. Thus Moonshine is related somehow to Lie theory.

McKay later found independent relationships with Lie theory [{\bf 89}],
[{\bf 15}], [{\bf 47}],
reminiscent of his famous A-D-E correspondence with finite subgroups of
SU$_2(\C)$. As mentioned earlier, $\M$ has two conjugacy classes of involutions.
Let $K$ be the smaller one, called `2A' in [{\bf 22}] (the alternative, class
`2B', has almost 100 million times more elements). The product of any two
elements of $K$ will lie in one of nine conjugacy classes: namely, 1A, 2A, 2B,
3A, 3C, 4A, 4B, 5A, 6A, corresponding respectively to elements of orders
1, 2, 2, 3, 3, 4, 4, 5, 6. It is surprising that,
for such a complicated group as $\M$, that list stops at only 6 --- we call
$\M$ a {\it 6-transposition group} for this reason (more on this in \S5.2).
The punchline: McKay noticed that those nine numbers are precisely the labels
of the affine $E_8$ diagram (see Figure 2). Thus we can attach a conjugacy
class of $\M$ to each vertex of the $\widehat{E}_8$ diagram.
An interpretation of the {\it
edges} in the $\widehat{E_8}$ diagram, in terms of $\M$, is unfortunately not
known.

\medskip\epsfysize=1.5in\centerline{ \epsffile{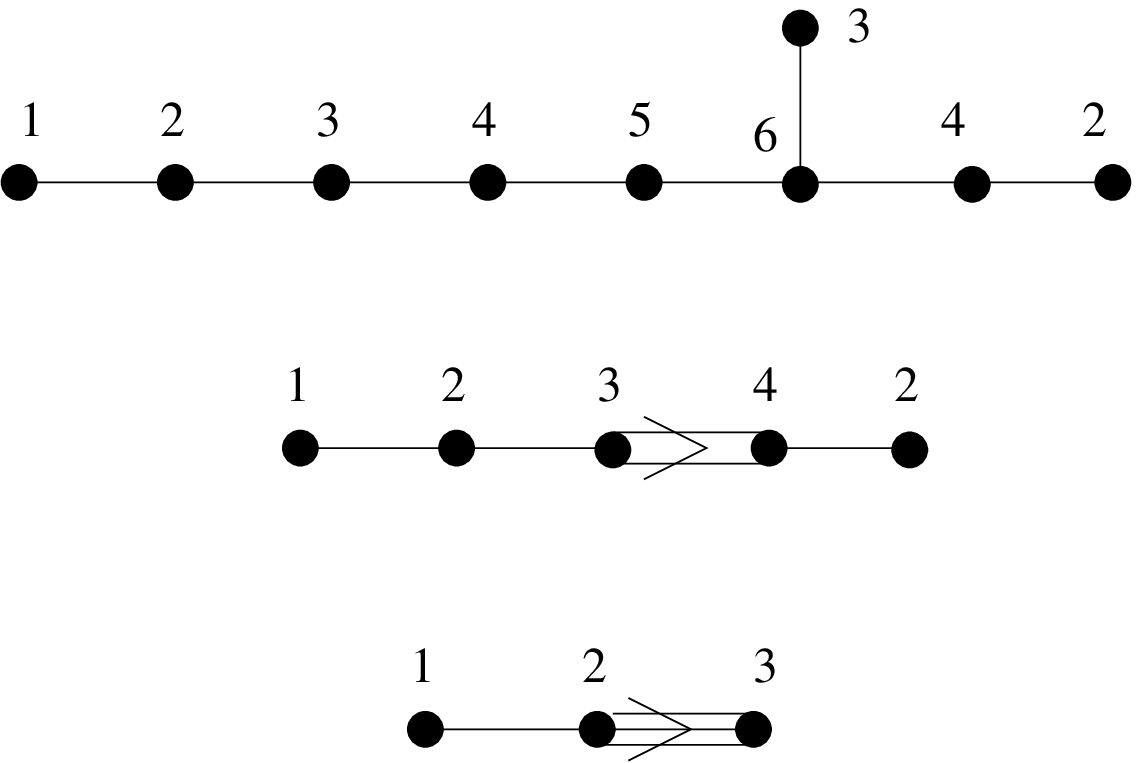}}\medskip
\centerline{{\bf Figure 2. The affine $E_8$, $F_4$, and $G_2$ diagrams with labels}}\medskip

We can't get the affine $E_7$ labels in a similar way, but McKay noticed that
an order {two} folding of affine $E_7$ gives the affine $F_4$ diagram, and we
can obtain its labels using the Baby Monster ${\Bbb B}$ (the second largest sporadic).
In particular, let $K$ now be the smallest conjugacy class of involutions in
${\Bbb B}$ (also labelled `2A' in [{\bf 22}]); the conjugacy classes in $KK$ have
orders 1, 2, 2, 3, 4 (${\Bbb B}$ is a 4-transposition group), and these are
the labels of $\widehat{F_4}$. Of course we'd prefer $\widehat{E_7}$ to
$\widehat{F_4}$, but perhaps that {\it two}-folding has something to do with
the fact that an order-{\it two} central extension of ${\Bbb B}$ is the
centraliser of an element $g\in\M$ of order {\it two}.

Now, the {\it triple}-folding of affine $E_6$ is affine $G_2$. The Monster has three
conjugacy classes of order {\it three}. The smallest of these (`3A') has a centraliser
which is a {\it triple} cover of the Fischer group $Fi_{24}'.2$. Taking the smallest
conjugacy class of involutions in $Fi_{24}'.2$, and multiplying it by itself,
gives conjugacy classes with orders 1, 2, 3 (hence $Fi_{24}'.2$ is a 3-transposition
group) --- and those not surprisingly are the labels of $\widehat{G_2}$!

Although we now understand (3.5) (see \S4.1) and have proven the fundamental
Conway--Norton conjecture (see \S\S4.2--4.4), McKay's $\widehat{E_8},\widehat{F_4},
\widehat{G_2}$ observations still have no explanation. In [{\bf 47}] these patterns
are extended, by relating various simple groups to the $\widehat{E_8}$ diagram
with deleted nodes.

\medskip{{\it 3.3. Replicable functions.}}
There are several other less important conjectures. One which played  an important
role in ultimately proving the main conjecture involves the {\it replication
formulae}. Conway--Norton want to think of the Hauptmoduls $T_g$ as
being intimately connected with $\M$; if so, then the group structure
of $\M$
should somehow directly relate different $T_g$. Considering the
power map $g\mapsto g^n$ leads to the following.

It was well-known classically that $j(\tau)$ has the property that
$j(p\tau)+j({\tau\over p})+j({\tau+1\over p})+\cdots+j({\tau+p-1\over p})$
is a polynomial in $j(\tau)$, for any prime $p$ ({\it proof:}
 it's a modular
function for SL$_2(\Z)$, and hence equals a rational function of $j(\tau)$; since its
only poles will be at the cusps, the denominator polynomial must be trivial).
Hence the same will hold for $J$. More generally, we get
$$\sum_{ad=n,0\le b<d}J({a\tau+b\over d})=Q_{n}(J(\tau))\ ,\eqno(3.6a)$$
where $Q_{n}$ is the unique polynomial for which
$Q_{n}(J(\tau))-q^{-n}$ has a $q$-expansion with only strictly positive
powers of  $q$. For example, $Q_2(x)=x^2-2a_1$ and $Q_3(x)=x^3-3a_1x-3a_2$,
where we write $J(\tau)=\sum_n a_nq^n$. The left side of (3.6a) is really a
 Hecke operator applied to $J$. These equations (3.6a) can be rewritten into 
recursions such as
$a_4=a_3+(a_1^2-a_1)/2$, or collected together into the remarkable expression
(originally due to Zagier)
$$p^{-1}\prod_{{m>0\atop n\in\Z}}(1-p^mq^n)^{a_{mn}}=J(z)-J(\tau)\ ,\eqno(3.6b)$$
where $p=e^{2\pi \i z}$.

Conway and Norton  conjectured [{\bf 24}] that these formulas
have an analogue for any McKay--Thompson series $T_g$. In particular,
(3.6a) becomes
$$\sum_{ad=n,0\le b<d}T_{g^a}({a\tau+b\over d})=Q_{n,g}(T_g(\tau))\ ,\eqno(3.7a)$$
where $Q_{n,g}$ plays the same role for $T_g$ that $Q_n$ played for $J$.
These are called  the replication formulae.  Again, these yield recursions like
$a_4(g)=a_2(g)+(a_1(g)^2-a_1(g^2))/2$, or can be collected into the expression
 $$p^{-1}\exp\bigl[-\sum_{k>0}\sum_{{m>0\atop
n\in\Z}}a_{mn}(g^k){p^{mk}q^{nk}\over k}\bigr]=T_g(z)-T_g(\tau)\ .\eqno(3.7b)$$
This looks a lot more complicated than (3.6b), but you can glimpse the Taylor
expansion of ln$(1-p^mq^n)$ there and in fact for $g=1$ (3.7b) reduces
to (3.6b). 

Axiomatising (3.7a) leads to Norton's notion of {\it replicable function}
[{\bf 96}], [{\bf 1}]. Write $f^{(1)}(\tau)=q^{-1}+\sum_{k=1}^\infty
b_k^{(1)}q^k$, and replacing each
$T_{g^a}$ in (3.7a) with $f^{(a)}$, use (3.7a) to recursively define each
$f^{(n)}$. If each $f^{(n)}$ has a $q$-expansion of the form $f^{(n)}(\tau)=
q^{-1}+\sum_{k=1}^\infty b_k^{(n)}q^k$ --- i.e.\ no fractional powers of $q$
arise --- then we call $f=f^{(1)}$ replicable. Equation (3.7a) says the McKay--Thompson
series are replicable, and [{\bf 30}] proved that the Hauptmodul of any genus 0 modular
group of moonshine-type is replicable, provided its coefficients are rational.
Conversely, Norton conjectured that any replicable function with rational
coefficients is either such a Hauptmodul, or one of the  `modular fictions'
$f(\tau)=q^{-1}$, $f(\tau)=q^{-1}\pm q$. This conjecture seems difficult and
is still open.
Incidentally, if the coefficients $b_k^{(1)}$ are irrational, then the definition (3.7a) of
replicability should be modified to include Galois automorphisms (see \S8 of
[{\bf 29}]). Replication in positive genus is discussed in [{\bf 109}].

Replication (3.7a) concerns the power map $g\mapsto g^n$ in $\M$. Can
Moonshine see more of the group structure of $\M$? We explored one step in
this direction in \S3.2, where  McKay modeled products of conjugacy
classes using Dynkin diagrams. A different idea is given in \S5.1.
It would be very desirable to find other direct connections between the
group operation in $\M$ and e.g.\ the McKay--Thompson series.

\medskip{{\it 3.4. The Leech lattice and Moonshine.}}
The Leech lattice $\Lambda=\Lambda_{24}$ is a 24-dimensional even self-dual lattice [{\bf 26}] which is to
lattices much as the $\M$-module $V$ of (3.2a) turns out to be for vertex operator
algebras (see \S4.2 below). $\Lambda$ has no vectors of odd norm, no norm-2 
vectors, and precisely 196560 norm-4
vectors --- a number remarkably close to the monstrous 196883. In fact its 
theta series $\Theta_{\Lambda}(\tau)
=\sum_{v\in\Lambda}q^{v\cdot v/2}$, when divided by $\eta(\tau)^{24}$, equals
$J(\tau)+24$. Is this another example of Moonshine?

Indeed it is. However we have:

\medskip{{\smcap Theorem 3.}} {\it Let $L\subset\R^n$ be any $n$-dimensional
positive-definite lattice whose norms ${\bf v}\cdot{\bf v}$ are all rational. Let
${\bf t}\in \R^n$ be any vector with finite order in $L$: i.e.\ $m{\bf t}
\in L$ for some nonzero $m\in\Z$. Then the theta series}
$$\Theta_{L+{\bf t}}(\tau):=\sum_{{\bf v}\in L}e^{\pi\i\tau\,({\bf v}
+{\bf t})^2}\ ,$$
{\it divided by $\eta(\tau)^n$, is a modular function  for some 
$\Gamma(N)$.}\medskip

See e.g.\ Theorem 20 of [{\bf 100}] for a proof of this classical result. If $L$ is in fact
an {\it even} lattice (i.e.\ all norms ${\bf v}\cdot{\bf v}$ lie in $2\Z$),
we can say more. Let $L^*:=\{{\bf x}\in\R^n\,|\,{\bf x}\cdot L\subseteq
\Z\}$ be the dual lattice. It contains $L$ with finite index; write
${\bf t}_i+L$, $i=1,\ldots,M$, for the finitely many cosets in $L^*/L$.
Define a column vector $\vec{\chi}_L(\tau)$ with $i$th component
$\Theta_{{\bf t}_i+L}(\tau)/\eta(\tau)^n$. Then $\vec{\chi}_L$ forms a {\it vector-valued
modular function} for SL$_2(\Z)$: for any $A=\left(\matrix{a&b\cr c&d}\right)
\in{\rm SL}_2(\Z)$,
$$\vec{\chi}_L({a\tau+b\over c\tau+d})=\rho(A)\,\vec{\chi}_L(\tau)\eqno(3.8)$$
 for some $M$-dimensional
unitary matrix representation $\rho$ of SL$_2(\Z)$. In particular, for the 
Leech lattice $L=\Lambda$, $M=1$ and we can quickly identify 
$\Theta_{\Lambda}(\tau)$ in terms of $J(\tau)$. Although the $196560\approx 196884$
coincidence is thus trivial
to explain, it will turn out to be a very instructive example of Moonshine.

The lattices are related to groups through their automorphism groups, which will
always be finite for positive-definite lattices. The automorphism group $Co_0:=
{\rm Aut}(\Lambda)$ 
of the Leech lattice has order about $8\times 10^{18}$, and
is the direct product of $C_2$ with Conway's sporadic group $Co_1$.
Several other sporadics are also involved in Aut($\Lambda)$.
To each automorphism $\alpha \in {\rm Aut}(\Lambda)$, let $\theta_\alpha$
denote the theta series of the sublattice of $\Lambda$ fixed by $\alpha$.
[{\bf 24}] also associate to each automorphism $\alpha$ a  certain function $\eta_\alpha(\tau)$
of the form $\prod_i \eta(a_i\tau)/\prod_j\eta(b_j\tau)$ built out of the
Dedekind eta. Both $\theta_\alpha$ and $\eta_\alpha$ are constant on each
conjugacy class in Aut($\Lambda)$, of which there are 202. [{\bf 24}] remarks that
the  ratio $\theta_\alpha/\eta_\alpha$ always seems to equal some McKay--Thompson
series $T_{g(\alpha)}$. See also [{\bf 87}].

It turns out that this observation isn't quite correct [{\bf 74}]. For each $\alpha\in{\rm
Aut}(\Lambda)$, the subgroup of SL$_2(\R)$ which fixes $\theta_\alpha/\eta_\alpha$
is indeed genus 0, but for exactly 15 conjugacy classes in Aut$(\Lambda)$,
$\theta_\alpha/\eta_\alpha$ is not the Hauptmodul. 

Similarly, one can ask this for the $E_8$ lattice, whose automorphism group
is the Weyl group (of order $\approx 7\times 10^8$) of the $E_8$ Lie group.
The automorphisms $\alpha$ of the
$E_8$ lattice which yield a Hauptmodul were classified in [{\bf 19}].

\bigskip\centerline{{\it 4. Proof of the Monstrous Moonshine conjectures}}\medskip

At first glance, the significance of the Moonshine conjectures seems very
unlikely: they constitute after all a finite set of very specialised coincidences. 
The whole point though
is to try to understand {\it why} such seemingly incomparable objects as the
Monster and the Hauptmoduls can be so related,
and to try to extend and apply this understanding to other contexts.
Establishing the truth (or falsity) of the conjectures was merely meant
as an aid to uncovering the why of Moonshine. Indeed, in  achieving this 
understanding,
important new algebraic structures were formulated. We will sketch this
theory below.

The main Conway--Norton conjecture was attacked almost immediately. Thompson
showed [{\bf 113}] (see also [{\bf 103}]) that if $g\mapsto a_n(g)$ is a character for all sufficiently small
$n$ (apparently $n\le 1300$ is sufficient), then
it will be for all $n$. He also showed that if certain congruence conditions
hold for a certain number of $a_n(g)$ (all with $n\le 100$), then all
$g\mapsto a_n(g)$ will be {\it virtual} characters (i.e.\ differences of
true characters of $\M$).
Atkin, Fong, and Smith (see [{\bf 110}] for details) used that to prove on a computer
that indeed all $a_n(g)$ were virtual characters (they didn't quite reach
$n=1300$ though). But their work doesn't say
anything more about the underlying (possibly virtual) representation $V$,
other than its existence.  Their work plays no role in the following.

We want to show that the McKay--Thompson series $T_g(\tau)$
 of (3.2c) equals the Hauptmodul $J_{G_g}(\tau)$ in (3.3). 
First, we need to construct the infinite-dimensional module $V$ of $\M$.
We discuss this, and the underlying theory of vertex operator algebras, in
\S4.2. Borcherds' strategy [{\bf 11}] was to
bring in Lie theory, by associating to the module $V$ a `Monster Lie algebra'.
This algebra, and the underlying theory of generalised Kac--Moody algebras,
is described in \S4.3. In the final subsection we go from the Monster Lie algebra
to the replication formulae, and conclude the proof.
We begin this section though by explaining the much simpler connection of
$E_8$ with $j^{{1\over 3}}$.

\medskip{{\it 4.1. $E_8$ and $j^{{1\over 3}}$.}}
An explanation for the relation between $E_8$ and $j^{{1\over 3}}$ was found
almost immediately, by Kac and Lepowsky [{\bf 65}], [{\bf 76}]: $j^{{1\over 3}}$ is the (normalised)
character of a representation of the affine Kac--Moody algebra $E^{(1)}_8$.
Given a finite-dimensional simple Lie algebra $\g$, the affine algebra
$\g^{(1)}$ is the infinite-dimensional Lie algebra consisting of all Laurent
polynomials $\sum_{n=-\infty}^\infty a_nt^n$ where $a_n\in\g$ and $t$ is an
indeterminate, together with
a central term and derivation $D=-L_0$ (see [{\bf 66}], [{\bf 70}]). Highest weight representations are defined in
the usual way. Thanks largely to the fact that the affine Weyl group is a semi-direct
product of the additive group $\Z^r$ ($r={\rm rank}(\g)$) with the finite
Weyl group, the characters of these representations
(especially the `integrable' highest weight ones, which are the direct analogue
of the finite-dimensional representations of $\g$) transform nicely with respect
to SL$_2(\Z)$. See e.g.\ Chapter 13 of [{\bf 70}] for details. This is probably the
single biggest reason Kac--Moody algebras are so well-known. 

\medskip{{\smcap Theorem 4.}} [{\bf 69}] {\it Let $\g$ be any finite-dimensional Lie algebra,
and $\g^{(1)}$ denote the corresponding affine algebra. Let $P_+^k$ denote the finitely
many `level $k$ integrable highest weight modules' $L_\lambda$ of $\g^{(1)}$.
The $\g^{(1)}$-module $L_\la$ has a natural $\Z$-grading $L_\la=\oplus_{n=0}^\infty(L_\la)_n$
into finite-dimensional $\g$-modules $(L_\la)_n$. Let $\chi_\la(\tau)=
q^{h_\la-c/24}\sum_{n=0}^\infty {\rm dim}((L_\la)_n) \,q^n$ be the corresponding normalised
character (for some appropriate choice of $h_\la-c/24\in\Q$). Then each
$\chi_\la$ is a holomorphic function in $\H$, and the vector $\vec{\chi}_k(\tau)$
with entry $\chi_\la(\tau)$ for each $\la\in P_+^k$ defines a vector-valued
modular function for SL$_2(\Z)$, as in (3.8), 
for some finite-dimensional unitary representation $\rho$ of SL$_2(\Z)$.}\medskip

In fact each character $\chi_\la$ will be a rational function in lattice theta
series, and so will be a modular function for some
$\Gamma(N)$. It turns out that there is only
one level 1 integrable highest-weight representation of $E_8^{(1)}$, and
its character equals $j^{{1\over 3}}$. The modularity of $j^{{1\over 3}}$ is
thus predicted
by Kac--Moody theory, and the fact that the coefficients are dimensions of
$E_8$ representations is automatic. We will see a simultaneous generalisation
of Theorems 3 and 4 next section.

We've already encountered the mysterious normalisation $q^{-1}$ of the McKay--Thompson
series, and the $q^{-{1\over 3}}$ of the $E_8^{(1)}$ character, and more
generally the $q^{h_\la-c/24}$ of $\chi_\la$. Many explanations have been
provided for this pervasive factor. For example, to [{\bf 2}] it is topological
in origin, and related to the Atiyah--Singer Index Theorem; a geometric
interpretation using determinant line bundles is due to Segal [{\bf 107}]. In 
quantum
physics it's called the {\it conformal anomaly} (a breakdown of manifest
conformal symmetry when the classical system is quantised), and is introduced
in regularisation as a vacuum energy. Probably the simplest instance of it is the prefactor in the
familiar definition $\eta(\tau)=q^{{1\over 24}}\prod_{n=1}^\infty(1-q^n)$
for the Dedekind eta: reading through classical proofs for its modularity
we find that `${1\over 24}$' here arises through the combination $\zeta(2)/
(2\pi)^2$; $\eta$ appears also in physics in the partition function of the bosonic string,
and that same `${1\over 24}$' arises there via regularisation as $-\zeta(-1)/2$.
The equivalence of these two expressions for ${1\over 24}$ comes from the 
functional equation of  the Riemann zeta.
This same zeta value appears famously in the central term of the Virasoro
algebra (4.3), and Bloch [{\bf 8}] found other zeta values appearing in other algebras
of differential operators, many of which have now been  interpreted and
generalised (starting with [{\bf 77}])
 within the vertex operator algebra framework.

Although a direct explanation for Monstrous Moonshine using affine algebras
has never been found (and certainly isn't expected), the theory of Kac--Moody algebras
influenced every stage
of the ultimate proof. For this reason we'll briefly sketch their theory.
A simple finite-dimensional Lie algebra $\g$ is  built out of the 3-dimensional
algebra sl$_2$, in a simple way; the Dynkin diagram of $\g$ encodes the
exact presentation. In the identical way, Kac--Moody algebras are also built out 
of copies of sl$_2$ --- the only difference is that the finite-dimensionality
constraint  (a positive-definiteness condition on the Cartan matrix) is lifted.
Their structure is completely analogous to that of the simple Lie algebras:
e.g.\ it has a grading by roots into finite-dimensional spaces; it has a
triangular decomposition (making Verma modules possible); it has an invariant
symmetric bilinear form. See e.g.\ [{\bf 66}], [{\bf 70}] for details. The
affine algebras are the class of Kac--Moody algebras
especially analogous to the finite-dimensional ones.

\medskip{{\it 4.2. The Moonshine module $V^\natural$.}}
A vital component of the Monstrous Moonshine conjectures came a few years
after [{\bf 24}]. In a deep work, Frenkel--Lepowsky--Meurman [{\bf 43}],
[{\bf 44}] constructed a graded infinite-dimensional
representation $V^{\natural}$ of $\M$ and conjectured (correctly) that it
is the representation $V$ in (3.2a). $V^\natural$ has a very rich algebraic
structure: it is in fact a {\it vertex operator algebra!}

A vertex operator algebra [{\bf 9}], [{\bf 44}], [{\bf 67}], [{\bf 39}],
[{\bf 79}] is an infinite-dimensional vector space
$V$ with infinitely many heavily constrained bilinear products $u*_nv$.
The name means `algebra of (generalised) vertex operators'; vertex operators
are formal differential operators which originally appeared in physics as quantum
fields describing the creation and propagation of physical strings (see 
\S6 below), and were
constructed later but independently by Lie theorists (starting with Lepowsky
and Wilson)
to realise affine Kac--Moody algebras as algebras of differential operators.
Because there were vertex operator constructions associated to lattices,
affine algebra modules, and string theory, and all of these have connections
to modular functions, it was natural to use vertex operators 
to try to construct the $\M$-module $V$ of (3.2a). 

The definition of vertex operator algebra (VOA) is too complicated to give in
detail here. A VOA is a graded infinite-dimensional vector
space ${\cal V}=\oplus_{n=0}^\infty {\cal V}_n$, where each ${\cal V}_n$ is finite-dimensional.
To simplify the discussion, we will limit ourselves in this paper to VOAs
with one-dimensional
${\cal V}_0$, which is typical of the examples relevant to Moonshine (and conformal
field theory).
To each vector $v\in {\cal V}$ we assign a {\it vertex operator} $Y(v,z)$, which is a 
formal power series $Y(v,z)=\sum_{m\in\Z}v_{(m)} z^{-m-1}$, with coefficients
$v_{(m)}\in{\rm End}({\cal V})$. The vertex operator is just the generating function
for the products: $u*_nv=u_{(n)}(v)$. These products respect the grading ---
in particular,
$${\cal V}_k*_n{\cal V}_\ell\subseteq {\cal V}_{k+\ell-n-1}\ .\eqno(4.1)$$
A key axiom, which collects together all the identities obeyed by the
products $u*_nv$, can be written as
$$(z-w)^M[Y(u,z),Y(v,w)]=0\qquad\forall u,v\in {\cal V}\ ,\eqno(4.2a)$$
for some integer $M$ (depending on $u,v$), where the bracket in (4.2a) means
the commutator $Y(u,z)\,Y(v,w)-Y(v,w)\,Y(u,z)$. This strange-looking formula
really says that each such commutator is a linear combination of Dirac deltas
and their derivatives, all centred at $z=w$ (see e.g.\ Corollary 2.2 in 
[{\bf 67}]). Equation (4.2a) implies more down-to-earth identities, such as 
$$(u_\ell v)_n=\sum_{i\ge 0}(-1)^i\left({\ell\atop i}\right)(u_{\ell-i}\circ
v_{n+i}-(-1)^\ell v_{\ell+n-i}\circ u_i)\ .\eqno(4.2b)$$

There are two distinguished elements in ${\cal V}$: the identity ${\bf 1}\in {\cal V}_0$
(so ${\cal V}_0=\C{\bf 1}$) and the {\it conformal vector} $\omega\in {\cal V}_2$. 
The identity obeys $Y({\bf 1},z)=id$,
i.e.\ ${\bf 1}_{(n)}v=\delta_{n,-1}v$. More interesting is the conformal vector:
writing $L_n=\omega_{(n+1)}$, the operators $L_n$ are required to form a representation
on ${\cal V}$ of the {\it Virasoro algebra}:
$$[L_m,L_n]=(m-n)L_{m+n}+\delta_{m,-n}{m^3-m\over 12}c \,id_{{\cal V}}\ ,\eqno(4.3)$$
for some number $c\in\R$ (an important numerical invariant of ${\cal V}$) called the 
{\it rank} or {\it central charge} of the VOA. In addition we
require $L_0u=nu$ whenever $u\in {\cal V}_n$, and $L_{-1}$ acts on ${\cal V}$ as a derivation.

The appearance here of the Virasoro algebra is fundamental. It is the 
unique nontrivial central extension of the (polynomial) vector fields on $S^1$,
which in turn is the Lie algebra associated to the group Diff($S^1$) of 
diffeomorphisms of the circle.

The notion of a VOA may seem very arbitrary, but as we'll mention in \S6 it
is the `chiral algebra' of a conformal field theory.
The simplest (and least interesting) special case of a VOA occurs when $M=0$
in (4.2a) --- i.e.\ when all vertex operators commute. Then it is not hard to see
that, for each choice of $z\ne 0$, ${\cal V}$ would be a commutative associative 
algebra with unit, whose product is given by $u*_zv:=Y(u,z)v$. A more honest
way to motivate VOAs has been suggested by Huang: binary trees can be used
to keep track of the brackets in nested products, e.g.\ $a((bc)d)$, and
e.g.\ Lie algebras can be easily formulated using this language [{\bf 58}]; in a monumental
work [{\bf 59}], Huang `two-dimensionalised' this Lie algebra formulation by replacing
binary trees with spheres with tubes, and showed that the result is equivalent
to a VOA.

A relation between VOAs and Lie algebras also exists at a more elementary level.
Placing $\ell=n=0$ in (4.2b), and evaluating it on the right by $w\in {\cal V}$,
gives 
$$(u_0v)_0w=u_0(v_0w)-v_0(u_0w)\ .\eqno(4.4)$$ 
Writing $[xy]$ for $x_0y$, this becomes
the Lie algebra Jacobi identity, at least when $[xy]=-[yx]$. There
are different ways to obtain from this a true Lie algebra. The simplest is 
that $V_1$ will be a Lie algebra for that bracket;
moreover, the number $\langle u,v\rangle$ defined by
$u_1v=\langle u,v\rangle{\bf 1}$ is an invariant bilinear form for this
Lie algebra (by (4.2b) with $\ell=0,n=1$). Equation (4.4) tells us each
${\cal V}_n$ is a ${\cal V}_1$-module, and in fact $e^u$ will be an automorphism of ${\cal V}$
for any $u\in {\cal V}_1$. In the most common examples, ${\cal V}_1$ will be
reductive (i.e.\ a direct sum of simple and abelian Lie algebras).

Modules $M$ of ${\cal V}$ can be defined in the obvious way [{\bf 42}],
[{\bf 79}] --- e.g.\ for each
$u\in {\cal V}$, each $u_{(n)}$ will be in End$(M)$. 
All (irreducible) modules $M$ come with a $\Z$-grading $M=\oplus_{n=0}^\infty
M_n$, where $u_kM_n\subseteq M_{n+\ell-k-1}$ for all $u\in {\cal V}_\ell$,
and $L_0x=(n+h)x$ for any $x\in M_n$, where $h$ is some number (the {\it conformal
weight}) depending only on $M$. The (normalised) character of $M$ is
$$\chi_M(\tau)=q^{-c/24}{\rm Tr}_M q^{L_0}=q^{h-c/24}\sum_{n=0}^\infty {\rm dim}
(M_n)\,q^{n}\ .$$

It takes a little effort to construct even the simplest examples of VOAs.
The best-behaved ones are called {\it rational VOAs} [{\bf 124}] (borrowing on terminology
from physics) and have only finitely many irreducible modules. A rational VOA
is associated to any even positive-definite lattice $L$, and their modules are
in one-to-one correspondence with the cosets $L^*/L$. Another
important example: to any affine nontwisted
Kac--Moody algebra and choice of positive integer $k$ (the `level'), the highest 
weight module $L_{k\Lambda_0}$ has a natural VOA structure, and its modules are
precisely the affine algebra modules $L_\la$ for each highest weight 
$\la\in P_+^k$. 

One of the deepest results in the theory of VOAs is due to Zhu:

\medskip{{\smcap Theorem 5.}} [{\bf 124}] {\it Let ${\cal V}$ be a rational VOA. Its 
characters $\chi_M(\tau)$ are holomorphic in ${\Bbb H}$, and the subspaces $M_n$ 
carry representations for Aut$({\cal V}$). Write $\vec{\chi}_{{\cal V}}(\tau)$ 
for the vector whose components are the characters $\chi_M(\tau)$ of irreducible
modules $M$. Then 
$\vec{\chi}_{{\cal V}}$ is a vector-valued modular function for SL$_2(\Z)$.}\medskip

It is believed that the characters $\chi_M(\tau)$ themselves will be modular
functions for some $\Gamma(N)$;  significant progress towards this
was made in [{\bf 5}] (see also [{\bf 71}]). 
The proof of Zhu's Theorem is much more difficult than that of Theorem 4, which
it generalises.

The automorphism group Aut$({\cal V})$ is by definition required to fix $\omega$, 
which is why it respects the grading of ${\cal V}$. Aut$({\cal V})$ is how group theory
impinges on VOA theory. Since the automorphism group Aut$({\cal V})$ of a VOA contains
$e^{{\cal V}_1}$ as a (normal) subgroup, Aut(${\cal V}$) can be finite only when ${\cal V}_1=0$.
Zhu's Theorem tells us that Moonshine (without the genus-0 aspect) will hold
between the group Aut$({\cal V})$ and the functions $\chi_M(\tau)$, for
 any rational VOA.

The most famous example of a VOA is the Moonshine module $V^\natural$ of
[{\bf 44}]. It is the orbifold of the Leech lattice VOA ${\cal V}_{\Lambda}$
by the $\pm 1$-symmetry of $\Lambda$,
which means it's the direct sum of two parts: an invariant part $V^\natural_+$
and a twisted part $V^\natural_-$ (more on this in \S5.1).  The orbifold serves
two purposes: it removes the constant
term `24'  from the graded dimension $J+24$ (hence the subspace $({\cal V}_{\Lambda})_1$) 
of ${\cal V}_{\Lambda}$; and it enhances the symmetry from the discrete part
of Aut$({\cal V}_{\Lambda})$,
which is an extension of $Co_0$ by $(C_2)^{24}$, to all of $\M$.

A major claim of [{\bf 44}]
was that $V^\natural$ is a `natural' structure (hence their notation).
Even so, this bipartite structure to $V^\natural$ complicates its study.
We have $V^\natural_0=\C{\bf 1}$, as usual, but the Lie algebra ${\cal V}_1=\{0\}$
is trivial. For any such VOA, the space ${\cal V}_2$ will be a commutative nonassociative
algebra with product $u\times v:=u_1v$ and identity ${1\over 2}\omega$.
For the Moonshine VOA $V^\natural$, this can be shown (with effort!) to
be the 196883-dimensional {Griess algebra} extended by an identity
element. From this, we find the automorphism group of $V^\natural$ to be
 the Monster $\M$. The only irreducible module for $V^\natural$ is itself
--- such a VOA is called {\it holomorphic}.
Together with Zhu's Theorem, this implies that its character, namely $J(\tau)$,
must be a modular function for SL$_2(\Z)$ (strictly speaking, we only get 
invariance up to a 1-dimensional character of SL$_2(\Z)$, but it is easy to
show that character must be identically 1). We'll see in \S5.1 how to obtain the other
McKay--Thompson series from $V^\natural$.

Conjecturally, there are  71 holomorphic VOAs with rank $c=24$ [{\bf 106}].
Much as the Leech lattice is the unique even self-dual positive-definite lattice of
dimension 24 containing no norm-2 vectors [{\bf 26}], the Moonshine module 
$V^\natural$ is {\it conjecturally} [{\bf 44}] the unique holomorphic VOA with
$c=24$ and with trivial ${\cal V}_1$.
Thus, just as the Leech lattice is the unique lattice with
theta series $\Theta_{\Lambda}$, so (conjecturally) is the Moonshine
module the unique holomorphic VOA with (normalised) graded dimension $J$. Proving this is one of the
most important (and difficult) challenges in the subject.

\medskip{{\it 4.3. The Monster Lie algebra ${\frak m}$}}. 
To show that all of the McKay--Thompson series $T_g$ are indeed Hauptmoduls,
Borcherds needed identities satisfied by their $q$-expansions. He obtained
these through a Lie algebra he associated to $V^\natural$. Before discussing 
it, let's briefly describe Borcherds' generalisation of Kac--Moody algebras
[{\bf 10}].

A Borcherds--Kac--Moody algebra differs from a Kac--Moody algebra
in that it is built up from Heisenberg algebras as well as sl$_2$, and
these subalgebras intertwine in more complicated ways. Nevertheless
much of
the theory for finite-dimensional simple Lie algebras continues to find
an analogue in this much more general setting (e.g.\ root-space decomposition,
Weyl group, character formula,...). This unexpected fact is the point
of Borcherds--Kac--Moody algebras. For reasons of space we avoid giving here the
fairly simple definition, but for this and much more see the review articles 
[{\bf 53}], [{\bf 63}], [{\bf 102}].

Their basic structure theorem is that of Kac--Moody algebras. In particular,
there is a grading by roots into finite-dimensional spaces (except that the
0-graded piece, corresponding to the Cartan subalgebra, may be infinite-dimensional).
They also have a triangularisable decomposition and an invariant symmetric
bilinear form. Indeed, these structural properties {\it characterise}
Borcherds--Kac--Moody algebras. In this sense
Borcherds--Kac--Moody algebras are the ultimate
generalisation
of simple Lie algebras, in  that any further generalisation would lose
some basic structural ingredient.

In short, Borcherds' algebras strongly resemble the Kac--Moody
ones and constitute a natural and nontrivial generalization. The
main  differences are that they can be generated by copies of
the 3-dimensional Heisenberg algebra as well as sl$_2$, and that there can be
imaginary simple roots. Borcherds introduced these algebras and developed their
theory in order to understand the Monster Lie algebra ${\frak m}$.

We want to construct ${\frak m}$ from the Moonshine module $V^{\natural}
=V^\natural_0\oplus V^\natural_1\oplus\cdots$. For later convenience, 
relabel its subspaces $V^i:=V^\natural_{i+1}$.
 Of course the obvious choice $V^{\natural}_1=V^0$ is 0-dimensional, so we must
modify $V^\natural$ first. Let  $II_{1,1}$ denote the even self-dual indefinite
lattice consisting of all pairs $(m,n)\in\Z^2$ with inner product
$(m,n)\cdot(m',n')=mn'+nm'$. Because it is indefinite, the usual construction
of a VOA from a lattice will fail here to produce a true VOA, but most properties
will be obtained. Call this near-VOA, ${\cal V}_{{1,1}}$.

The Monster Lie algebra ${\frak m}$ is a Lie algebra associated to the near-VOA $V^\natural
\otimes  {\cal V}_{{1,1}}$ --- see [{\bf 11}] for the details. 
${\frak m}$ inherits a $II_{1,1}$-grading
from ${\cal V}_{{1,1}}$, and this is its root space decomposition: the 
$(m,n)$ root space
is isomorphic (as a vector space) to $V^{mn}$, if $(m,n)\ne (0,0)$;
the $(0,0)$ piece is isomorphic to $\R^2$. Structurally, 
the Monster Lie algebra has a decomposition ${\frak m}=u^+\oplus{\rm gl}_2\oplus
u^-$ into a sum of Lie subalgebras, where $u^{\pm}$ are free Lie algebras
(see e.g.\ [{\bf 64}]). It inherits the action of $\M$ from $V^\natural$.

This construction of ${\frak m}$ may seem indirect; an alternate approach, 
anticipated in [{\bf 11}] and [{\bf 12}], uses
{\it Moonshine cohomology} [{\bf 81}] ---  a functor, inspired by BRST
cohomology in conformal field theory, assigning
to certain $c=2$ near-VOAs some Lie algebra
carrying an action of $\M$. To ${\cal V}_{1,1}$ this functor associates ${\frak m}$.

\medskip{{\it 4.4. Denominator identities and modular equations.}}
It was discovered early on that  the
Hauptmoduls all obey the replication formulae, and that anything obeying
those formulae will be determined by their first few coefficients. The idea
then is to show that the McKay--Thompson series $T_g$ of (3.2c) also are
replicable. Borcherds did this using Lie algebra denominator identities [{\bf 11}].

Finite-dimensional simple Lie algebras ${\frak g}$ possess a very useful formula for their characters,
due to Weyl: the (formal) character $\chi_\la$ of a module $L_\la$ equals
$$\chi_\la:=\sum_{\mu}{\rm dim}(L_\la(\mu))\,e^{\mu
}=e^{-\rho}\,{\sum_{w\in W}\eps(w)\, e^{w(\la+\rho)}\over 
\prod_{\alpha\in\Delta_+} (1-e^{-\alpha})}\ ,\eqno(4.5)$$
where $W$ is the Weyl group, $\Delta_+$ the positive roots, $\eps(w)={\rm det}(w)$
is a sign, and where
$\oplus_\mu L_\la(\mu)$ is the weight-space decomposition of $L_\la$.
As the weights $\mu$ by definition lie in the dual ${\frak h}^*$ of the Cartan
subalgebra of ${\frak g}$, the character $\chi_\la$ can be regarded as a
complex-valued function on the space ${\frak h}\cong{\Bbb C}^r$ ($r={\rm rank}
({\frak g})$).

 Consider
the trivial representation: i.e.\ $x\mapsto 0$ for all $x\in \g$. Its
character $\chi_0$ will be identically 1. Thus the
character formula (4.5) tells us that a certain alternating  sum over a Weyl group, equals
a certain product over positive roots. These formulas, called {\it
denominator identities}, are nontrivial even in this finite-dimensional case.

In a famous paper [{\bf 83}], Macdonald generalised the denominator
identity for (4.5), to infinite sum/product identities, 
corresponding to the extended Dynkin diagrams. The simplest one was known
classically as the Jacobi triple product identity:
$$\sum_{n=-\infty}^\infty(-1)^nx^{n^2}y^n=
\prod_{m=1}^\infty(1-x^{2m})(1-x^{2m-1}y)(1-x^{2m-1}y^{-1})\ .\eqno(4.6)$$
Macdonald's identities were later reinterpreted, by Kac and Moody, as denominator
identities for the affine  algebras. For example, we
 now know (4.6) to be the denominator identity for the algebra $A_1^{(1)}$.

In particular,
the same formula (4.5) holds for Kac--Moody algebras, except that the
sum and product are now infinite, the positive roots now come with
multiplicities, and the characters are usually normalised by a prefactor 
$q^{h_\la-c/24}$. 
The variable $\tau$ in Theorem 4 is one of the coordinates in the 
Cartan subalgebra ${\Bbb C}^{r+2}$ of the affine algebra (see e.g.\ equation
(13.2.4) of [{\bf 66}]). In that theorem we dropped the
remaining variable dependence of the $\chi_\lambda$ for readability, although
those additional coordinates serve the important role of guaranteeing linear
independence of the characters, and of giving us an action of SL$_2(\Z)$
rather than merely PSL$_2(\Z)$.

Because a Borcherds--Kac--Moody algebra $\g$ is triangularisable, highest weight
$\g$-modules can be defined in the usual way from Verma modules.
The character formula becomes
$$\chi_{\la}=e^{-\rho}{\sum_{w\in W}\epsilon(w)\,w(e^{\la+\rho}S_\lambda)\over
\prod_{\alpha\in\Delta_+}(1-e^{-\alpha})^{{\rm mult}\,\alpha}}\ ,\eqno(4.7)$$
where  $S_\lambda$ is a correction factor due to imaginary simple roots.

The corresponding denominator identity of the Monster Lie algebra ${\frak m}$ can be
computed, and is given in (3.6b). Its Weyl group is $C_2$ and sends the
$(m,n$)-root space to $(n,m$); the $(m,n)$ root has multiplicity given by
coefficient $a_{mn}$ of $J$; for each $n>0$ we have an imaginary simple root
$(1,n)$ with multiplicity $a_{n}$. Because of a cohomological interpretation
of all denominator identities, (3.6b) can be `twisted' by each $g\in\M$,
and this gives (3.7b).
These formulas are equivalent to
the replication formula (3.7a) conjectured in \S 3.3.

Identities equivalent to (3.7b) were obtained by more elementary means --- i.e.\
methods requiring less of the theory of Borcherds--Kac--Moody algebras --- in
[{\bf 64}] and [{\bf 68}], permitting a simplification of Borcherds' proof
at this stage.

Now, it turns out that if we verify for each conjugacy class $K_g$ of $\M$
that the first, second, third, fourth and sixth coefficients of the
McKay--Thompson series $T_g$ and the corresponding Hauptmodul $J_{G_g}$ agree,
then indeed $T_g=J_{G_g}$.
That is precisely what Borcherds then did: he compared finitely many
coefficients, and as they all equalled what they should, this concluded
the proof [{\bf 11}] of Monstrous Moonshine!\medskip

However, this case-by-case verification
 occurred at the critical point where the McKay--Thompson series
were being compared directly to the Hauptmoduls, and so provides little insight
into why the $T_g$ are genus 0. Fortunately  a more conceptual
explanation of their equality has since been found.

A function $f$ obeying the replication formulae (3.7a) will also obey 
{\it modular equations} --- i.e.\ a 2-variable polynomial
identity satisfied by $f(x)$ and $f(nx)$. The simplest examples come from
the exponential and cosine functions: note that for any $n>0$,
$\exp(nx)=(\exp(x))^n$ and $\cos(nx)=T_n(\cos(x))$ where $T_n$ is a
Tchebychev polynomial. It was known classically that $j$ (hence $J$) satisfied
a modular equation for any $n$: e.g.\ put $X=J(\tau)$ and 
$Y=J(2\tau)$, then
$$\eqalignno{(X^2-Y)(Y^2-X)=&\,393768\,(X^2+Y^2)+42987520\,XY+40491318744\,(X+Y)
&\cr&\,-120981708338256\ .&\cr}$$

The only functions $f(\tau)=q^{-1}+a_1q+\cdots$
which obey modular equations for all $n$, are $J(\tau)$ and the `modular
fictions' $q^{-1}$ and $q^{-1}\pm q$ (which are essentially exp, cos, and
sin) [{\bf 72}]. More generally, we have:

\medskip{\smcap Theorem 6.} [{\bf 29}] {\it A function $B(\tau)
=q^{-1}+\sum_{n=1}^\infty b_nq^n$ which obeys a modular equation for all $n\equiv 1$
(mod $N$), will either be of the form $B(\tau)=q^{-1}+b_1q$, or will be a
Hauptmodul for a modular group of moonshine-type.} \medskip

The converse is also true [{\bf 29}]. The denominator identity argument tells us each
$T_g$ obeys a modular equation for each $n\equiv 1$ modulo the order
of $g$, so Theorem 6 then concludes the proof of Monstrous Moonshine.

The computer searches in [{\bf 20}] suggest that the hypothesis of
Theorem 6 may be considerably weakened, perhaps all the way down to the 
existence of modular equations for any two distinct primes.

\bigskip\centerline{{\it 5. Further developments}}

\medskip{{\it 5.1. Orbifolds.}} About a third of the McKay--Thompson series $T_g$ will
have some negative coefficients. In \S5.4 we'll see Borcherds interpret them
as dimensions of superspaces (which come with signs).
In an important announcement [{\bf 97}], on par with [{\bf 24}], Norton proposed
that, although $T_g(-1/\tau)$ will not usually be another McKay--Thompson series, it will
always have nonnegative integer $q$-coefficients, and these can be interpreted as
ordinary dimensions. In the process, he extended the $g\mapsto T_g$ assignment
to commuting pairs
$(g,h)\in \M\times\M$. 

In particular, to each such pair we have a
function $N(g,h;\tau)$, which we will call a {\it Norton series}, such that
$$N(g^ah^c,g^bh^d;\tau)=\alpha\,N(g,h;{a\tau+b\over c\tau+d})\qquad
\forall\left(\matrix{a&b\cr c&d}\right)\in {\rm SL}_2(\Z)\ ,\eqno(5.1)$$ 
for some root of unity
$\alpha$ (of order dividing 24, and depending on $g,h,a,b,c,d$).
The Norton series $N(g,h;\tau)$ is either constant, or generates the modular 
functions for a genus-0 subgroup of SL$_2(\Z)$ containing some $\Gamma(N)$ 
(but otherwise not necessarily of moonshine-type). Constant $N(g,h;\tau)$
arise when all elements of the form $g^ah^b$ (for gcd($a,b)=1$) are 
`non-Fricke' (an element
$g\in\M$ is called {\it Fricke} if the group $G_g$ contains an element
sending 0 to $\i\infty$ --- the identity 1 is Fricke, as are 120 of the 171 
$G_g$). Each $N(g,h;\tau)$ has a $q^{{1\over N}}$-expansion for that $N$; the 
coefficients of this expansion
 are characters evaluated at $h$ of some central extension of the centralizer 
$C_{\M}(g)$. Simultaneous conjugation of $g,h$ leaves the Norton series 
unchanged: $N(aga^{-1},aha^{-1};\tau)=N(g,h;\tau)$. 

For example, when $\langle g,h\rangle\cong C_2\times C_2$ and $g,h,gh$ are
all in class 2A, then $N(g,h;\tau)=\sqrt{J(\tau)-984}$.
The McKay--Thompson series are
recovered by the $g=1$ specialisation: $N(1,h;\tau)=T_h(\tau)$. 
This action (5.1) of SL$_2(\Z)$ is related to its natural action
on the fundamental group $\Z^2$ of the torus, as we'll see in \S6, as well as
a natural action of the braid group, as we'll see next subsection. 
Norton arrived at his conjecture empirically, by studying the data of Queen
(see  \S5.3).

The basic tool we have for approaching Moonshine conjectures is the theory
of VOAs, so we need to understand Norton's suggestion from that point of
view. For reasons of space, we'll limit this discussion to $V^\natural$, but 
it generalises. Given any automorphism $g\in{\rm Aut}(V^\natural)$, we can 
define $g$-twisted modules in the obvious way [{\bf 36}]. Then for each $g\in\M$,
there is a unique $g$-twisted module, call it $V^\natural(g)$, for $V^\natural$ 
--- this statement generalises
the holomorphicity of $V^\natural$ mentioned in \S4.2. 
More generally, given any automorphism $h\in{\rm Aut}(V^\natural)$ commuting
with $g$, $h$ will yield an automorphism of $V^\natural(g)$, so we can 
perform Thompson's twist (3.2c) and write 
$$q^{-c/24}{\rm Tr}_{V^\natural(g)}h\, q^{L_0}=:{\cal Z}(g,h;\tau)\ .\eqno(5.2)$$
These ${\cal Z}(g,h)$'s can be thought of as the building blocks of the graded
dimensions of various eigenspaces in $V^\natural(g)$: e.g.\ if $h$ has
order $m$, then the subspace of $V^\natural(g)$ fixed by automorphism $h$ 
will have graded dimension $m^{-1}\sum_{i=1}^m{\cal Z}(g,h^i)$. 
In the case of the Monster considered here, we have ${\cal Z}(g,h)=N(g,h)$.

The important paper [{\bf 36}] proves that, whenever the subgroup $\langle g,h\rangle$
generated by $g$ and $h$ is cyclic, then $N(g,h)$ will be a 
Hauptmodul satisfying (5.1). One way this will happen of course is
 whenever the orders of 
$g$ and $h$ are coprime. Extending [{\bf 36}] to all commuting 
pairs $g,h$ is one of the most pressing tasks in Moonshine.

This orbifold construction is the same as was used to construct $V^\natural$
from ${\cal V}_{\Lambda}$: $V^\natural$ is the sum of the `$\iota$'-invariant
subspace $V^\natural_+$ of ${\cal V}_{\Lambda}$ with the `$\iota$'-invariant subspace 
$V^\natural_-$ of the unique `$-1$'-twisted module for ${\cal V}_{\Lambda}$, 
where $\iota\in{\rm Aut}(\Lambda)$ is some involution. The
graded dimensions of $V^\natural_{\pm}$ are $2^{-1}({\cal Z}(\pm 1,1)+{\cal Z}
(\pm 1,\iota))$, respectively, and these sum to $J$.

The orbifold construction is also involved in an interesting reformulation
of the Hauptmodul property, due to Tuite [{\bf 116}]. Assume the uniqueness
conjecture: $V^\natural$ is the only VOA with graded dimension $J$. He argues
from this that, for each $g\in\M$, $T_g$ will be a Hauptmodul iff the only
orbifolds of $V^\natural$ are ${\cal V}_{\Lambda}$ and $V^\natural$ itself. In e.g.\ [{\bf
62}], this analysis is extended to some of Norton's $N(g,h)$'s,
where the subgroup $\langle g,h\rangle$ is not cyclic (thus going beyond [{\bf 36}]),
although again assuming the uniqueness conjecture. 

\medskip{{\it 5.2. Why the Monster?}}
That $\M$ is associated with {\it modular functions} can be explained
by it being the automorphism group of the Moonshine VOA $V^\natural$.
 But what is so special about this group $\M$ that
these modular functions $T_g$ and $N(g,h)$ should be Hauptmoduls?
This is still open. One approach is due to
Norton, and was first (rather cryptically) stated in [{\bf 97}]: the Monster is probably the
largest (in a sense) group with the 6-transposition property. Recall from
\S3.2 that a $k$-transposition group $G$ is one generated by a conjugacy class
$K$ of involutions, where the product $gh$ of any two elements of $K$ has order
$\le k$. For example, taking $K$ to be the transpositions in the symmetric
group $G=S_n$, we find that $S_n$ is 3-transposition.

A transitive
action of $\Gamma:={\rm PSL}_2(\Z)$ on a finite set $X$ with one distinguished
point $x_0\in X$,
is equivalent to specifying a finite index subgroup $\Gamma_0$ of $\Gamma$.
In particular, $\Gamma_0$ is the stabiliser $\{g\in\Gamma\,|\,g.x_0=x_0\}$
of $x_0$, $X$ can be identified with the cosets $\Gamma_0\backslash \Gamma$,
and $x_0$ with the coset $\Gamma_0$. (If we avoid specifying $x_0$, then
$\Gamma_0$ will be identified only up to conjugation.)

To such an action, we can associate an interesting triangulation of the closed surface
$\Gamma_0\backslash\overline{\H}$, called a (modular) {\it quilt}. The definition,
originally due to Norton and further developed by Parker, Conway, and Hsu,
is somewhat involved and will be avoided here (but see especially Chapter 3
of [{\bf 57}]). It is so-named because there is a polygonal `patch' covering every
cusp of $\Gamma_0\backslash\H$, and the closed surface is  formed by sewing together
the patches along their edges (`seams'). There are a total of $2n$ triangles
and $n$ seams in the triangulation, where $n$ is the index $\|\Gamma_0\backslash
\Gamma\|=\|X\|$.  The boundary of each patch has an even number of edges,
namely the double of the corresponding cusp width. The familiar formula
$$\gamma={n\over 12}-{n_2\over 4}-{n_3\over 3}-{n_\infty\over 2}+1$$
for the genus $\gamma$ of $\Gamma_0\backslash\H$ in terms of the index $n$ and the
numbers $n_i$ of $\Gamma_0$-orbits of fixed-points of order $i$, can be
interpreted in terms of the data of the quilt (see (6.2.3) of [{\bf 57}]),
and we find in particular that if every patch of the quilt has at most 6
sides, then the genus will be 0 or 1, and genus 1 only exceptionally.

In particular, we're interested in one class of these $\Gamma$-actions (actually
an SL$_2(\Z)$-action, but this doesn't matter). Recall that the braid group
$B_3$ has presentation
$$\langle \sigma_1,\sigma_2\,|\,\sigma_1\sigma_2\sigma_1=\sigma_2\sigma_1\sigma_2
\rangle\ ,\eqno(5.3a)$$
and centre ${Z}=\langle (\sigma_1\si_2\si_1)^2\rangle$ [{\bf 7}]. It
is related to the modular group by
$$B_3/{Z}\cong {\rm PSL}_2(\Z)\ ,\quad B_3/\langle (\sigma_1\si_2\si_1)^4
\rangle\cong {\rm SL}_2(\Z)\ .\eqno(5.3b)$$

Fix a finite group $G$ (we're most interested in the choice $G=\M$). 
We can define a right action of $B_3$ on triples $(g_1,g_2,g_3)\in G^3$ by
$$(g_1,g_2,g_3)\si_1=(g_1g_2g_1^{-1},g_1,g_3)\ ,\qquad(g_1,g_2,g_3)\si_2=
(g_1,g_2g_3 g_2^{-1},g_2)\ .\eqno(5.4a)$$
We will be interested in this action on the subset of $G^3$ where all
$g_i\in G$ are involutions. The action (5.4a) is equivalent to a reduced version, where
we replace $(g_1,g_2,g_3)$ with $(g_1g_2,g_2g_3)\in G^2$. Then (5.4a) becomes
$$(g,h)\si_1=(g,gh)\ ,\qquad(g,h)\si_2=(gh^{-1},h)\ .\eqno(5.4b)$$
These $B_3$ actions come from specialisations of the Burau and reduced Burau 
representations [{\bf 7}], respectively, and generalise to actions of
$B_n$ on $G^n$ and $G^{n-1}$. We can get an action of SL$_2(\Z)$ from the
$B_3$ action (5.4b) in two ways: either

\smallskip\item{(i)} by restricting to commuting pairs $g,h$; or

\smallskip\item{(ii)} by identifying
each pair $(g,h)$ with all its conjugates $(aga^{-1}, aha^{-1})$.

\smallskip\noindent Norton's
SL$_2(\Z)$ action of \S5.1 arises from the $B_3$ action (5.4b), when we
perform both (i) and (ii). 

The quilt picture was designed for this SL$_2(\Z)$ action.
The point of this construction is that the number of sides in each patch
is determined by the orders of the corresponding elements $g,h$. 
If $G$ is say a
6-transposition group (such as the Monster), and we take the involutions
$g_i$ from 2A, then each patch will have $\le 6$ sides, and the corresponding
genus will be 0 (usually) or 1 (exceptionally if at all). In this way we can relate
the Monster with a genus-0 property.

Based on the actions (5.4), Norton anticipates some analogue of Moonshine
valid for noncommuting pairs. CFT considerations (`higher genus orbifolds')
alluded to in \S6 suggest that  more natural should be 
e.g.\ quadruples $(g,g',h,h')\in\M^4$ obeying $ghg^{-1}h^{-1}=h'g'h'{}^{-1}
g'{}^{-1}$.

\smallskip
An interesting question is, how much does Monstrous Moonshine determine the
Monster? How much of $\M$'s structure can be deduced from e.g.\ McKay's
$\widehat{E}_8$ Dynkin diagram observation, and/or the (complete) replicability
of the $T_g$, and/or Norton's conjectures in \S5.1, and/or Modular Moonshine
in \S5.4 below? A small start toward this is taken in
 [{\bf 99}], where some control on the subgroups of $\M$ isomorphic to $C_p\times
C_p$ ($p$ prime) was obtained, using only the properties of the series $N(g,h)$.
For related work, see Chapter 8 of [{\bf 57}].

\medskip{{\it 5.3. Other finite groups.}}
It is natural to ask about Moonshine for other groups.
Indeed, the Hauptmodul for $\Gamma_0(2)+$ looks like
$$q^{-1}+4372q+96256q^2+12\,40002q^3+\cdots\eqno(5.5a)$$
and we find the relations 
$$4372=4371+1\ ,\qquad 96256=96255+1\ ,\qquad 12\,40002=11\,39374+4371+2\cdot 1
\ ,\eqno(5.5b)$$
where 1, 4371, 96255, and $11\,39374$ are all dimensions of irreducible
representations of the Baby Monster ${\Bbb B}$. Thus we find Moonshine
for ${\Bbb B}$! We will return to this example shortly.

Of course any subgroup of $\M$ automatically inherits Moonshine by
restriction, but obviously this isn't interesting. 
 Most constructions of the Leech lattice start with Mathieu's sporadic
 $M_{24}$ (see e.g.\ Chapters 10 and 11 of [{\bf 26}]), and most
constructions of the Monster involve the Leech lattice. Thus we are led 
to the following natural hierarchy of (most) sporadics:
\smallskip \item{(1)} $M_{24}$ (from which we can get $M_{11}$, $M_{12}$,
$M_{22}$, $M_{23}$); which leads to

\item{(2)} $Co_0=C_2\times Co_1$ (from which we get $HJ$,  $HS$,
$McL$, $Suz$, $Co_3$, $Co_2$); which leads to

\item{(3)} $\M$ (from which we get $He$, $Fi_{22}$, $Fi_{23}$, $Fi_{24}'$,
$HN$, $Th$, ${\Bbb B}$).\smallskip

It can thus be argued that we could approach problems in Monstrous Moonshine,
by first addressing in order $M_{24}$ and $Co_1$, which should be much simpler. Indeed, 
the full VOA orbifold theory --- i.e.\ the complete analogue of \S5.1 ---
for $M_{24}$ has been established in [{\bf 38}] (the relevant series
${\cal Z}(g,h)$ had already been constructed in [{\bf 88}]). 

Largely by trial and error, Queen [{\bf 101}] established Moonshine for the
following groups (all essentially centralisers of elements of $\M$):
$Co_0$, $Th$, $3.2.Suz$, $2.HJ$, $HN$, $2.A_7$, $He$, $M_{12}$ (by e.g.\
`$2.HJ$' we mean $C_2$ is normal and $HJ$ is the quotient $2.HJ/C_2$). 
In particular,
to each element $g$ of these groups, there corresponds a series $Q_g(\tau)=
q^{-1}+\sum_{n=0}^\infty a_n(g)q^n$, which is a Hauptmodul for some
modular group of moonshine-type, and where each $g\mapsto a_n(g)$ is a 
virtual character. For $Th$, 
$HN$, $He$ and $M_{12}$ it is a proper character. Other differences with
Monstrous Moonshine are that there can be a preferred nonzero value for the
constant term $a_0$, and that although $\Gamma_0(N)$ will be a subgroup of
the fixing group, it won't necessarily be normal. We will return to these
results next section, where we will see that many seem to come out of
the Moonshine for $\M$. About half of Queen's Hauptmoduls $Q_g$ for $Co_0$ do
not arise as a McKay--Thompson series for $\M$. Norton's conjectures in \S5.1
are a reinterpretation and extension of Queen's work.

Queen never reached ${\Bbb B}$ because of its size. However,  the Moonshine
(5.5) for ${\Bbb B}$ falls into her and Norton's scheme because (5.5a)
is the McKay--Thompson series associated to class 2A of $\M$, and the centraliser
of an element in 2A is a double cover of ${\Bbb B}$.

There can't be a VOA $V=\oplus_nV_n$ with graded dimension (5.5a) and automorphisms
in ${\Bbb B}$, because e.g.\ the ${\Bbb B}$-module $V_3$ doesn't contain $V_2$ as 
a submodule. However, H\"ohn deepened the analogy between $\M$ and ${\Bbb B}$ 
by constructing a vertex operator superalgebra $V{\Bbb B}^\natural$ of rank
$c=23.5$, called
the {\it shorter Moonshine module}, closely related to $V^\natural$ (see e.g.\
[{\bf 56}]). Its automorphism group is 
$C_2\times {\Bbb B}$. Just as $\M$ is the automorphism group of the Griess
algebra $V^\natural_2$, so is ${\Bbb B}$ the automorphism group of the algebra
$(V{\Bbb B}^\natural)_2$. Just as $V^\natural$ is associated to the Leech lattice
$\Lambda$, so is $V{\Bbb B}^\natural$ associated to the shorter Leech lattice $O_{23}$,
the unique 23-dimensional positive-definite self-dual lattice with no vectors 
of length 2 or 1 (see e.g.\ Chapter 6 of [{\bf 26}]). The automorphism group of
$O_{23}$ is $C_2\times Co_2$.

There has been no interesting Moonshine rumoured for the remaining six 
sporadics (the {\it pariahs} $J_1$, $J_3$, $Ru$, $ON$, $Ly$, $J_4$).
There will be some sort of Moonshine for any group which is an
automorphism group of a vertex operator algebra (so this means 
any finite group [{\bf 37}]!).  Many finite groups of
Lie type should arise as automorphism groups of VOAs associated to
affine algebras except defined over finite fields.
But apparently all known examples of genus-0 Moonshine are limited to the
groups involved with $\M$.

\medskip{{\it 5.4. Modular Moonshine}}.
Consider an element $g\in\M$. We expect from  [{\bf 101}], [{\bf 97}], [{\bf 36}]
that there is a Moonshine for the centraliser $C_{\M}(g)$ of $g$ in $\M$,
governed by the $g$-twisted module $V^\natural(g)$. Unfortunately, $V^\natural
(g)$ is not usually itself a VOA, so the analogy with $\M$ is not perfect. 
Ryba found it interesting that, for $g\in\M$ of prime order $p$, Norton's 
series $N(g,h)$ is a McKay--Thompson series (and has all the associated nice 
properties) whenever $h$ is $p$-{\it regular}
(i.e.\ $h$ has order coprime to $p$). This special behaviour of $p$-regular 
elements suggested to him to look at modular representations.

The basics of modular representations and {Brauer characters} are discussed in 
sufficient detail in Chapter 2 of [{\bf 31}]. A {\it modular representation} $\rho$
of a group $G$
is a representation defined over a field  of positive characteristic $p$ dividing
the order $|G|$ of $G$. Such representations possess many special (i.e.\
unpleasant) features. For one
thing, they are no longer completely reducible (so the role of irreducible 
modules as direct summands will be replaced with their role as composition
factors). For another, the usual notion
of character (the trace of representation matrices) loses its usefulness and is
replaced by the more subtle {\it Brauer character} $\beta(\rho)$: a complex-valued class function
on $\M$ which is only well-defined on the $p$-regular elements of $G$.

\medskip{\smcap Theorem 7.} [{\bf 105}], [{\bf 17}], [{\bf 13}] {\it Let
$g\in \M$ be any element of prime order $p$, for any $p$ dividing $|\M|$. 
Then there is a vertex operator
superalgebra ${}^g{\cal V}=\oplus_{n\in\Z}{}^g{\cal V}_n$ defined over the 
finite field ${\Bbb F}_p$ and acted on by the centraliser $C_{\M}(g)$. If $h\in C_\M(g)$
is $p$-regular, then the graded Brauer character 
$$R(g,h;\tau):=q^{-1}\sum_{n\in\Z}\beta({}^g{\cal V}_n)(h)\, q^n$$
equals the McKay--Thompson series $T_{gh}(\tau)$. Moreover, for $g$ belonging to
any conjugacy class in $\M$ except 2B, 3B, 5B, 7B, or 13B,  this is in fact
an ordinary VOA (i.e.\ the `odd' part vanishes),
while in the remaining cases the graded Brauer characters of both the
odd and even parts can separately be expressed using McKay--Thompson series.}\medskip

By a vertex operator superalgebra, we mean there is a $\Z_2$-grading into
even and odd subspaces, and for $u,v$ both odd the commutator in (4.2a) is
replaced by an anticommutator. In the proof, the superspaces arise as
cohomology groups, which naturally form an alternating sum.
The centralisers $C_\M(g)$ in the Theorem are quite nice: e.g.\ for $g$ in
classes 2A, 2B, 3A, 3B, 3C, 5A, 5B, 7A, 11A, respectively, these involve
 the sporadic
groups ${\Bbb B}$, $Co_1$, $Fi_{24}'$, $Suz$, $Th$, $HN$, $HJ$, $He$, and $M_{12}$.
The proof for $p=2$ is not complete at the present time.
The conjectures in [{\bf 105}] concerning modular
analogues of the Griess algebra for several sporadics follow from Theorem 7.

Can these modular ${}^g{\cal V}$s be interpreted as a reduction mod $p$ of 
(super)algebras in characteristic 0? Also, what about elements $g$ of composite
order?

\medskip{{\smcap Conjecture 8.}} [{\bf 13}] {\it Choose any $g\in\M$ and let
$n$ denote its order. Then there is a ${1\over n}\Z$-graded superspace
${}^g\widehat{{\cal V}}=\oplus_{i\in {1\over n}\Z}{}^g\widehat{{\cal V}}_i$
over the ring of cyclotomic integers $\Z[e^{2\pi\i/n}]$. It is often (but 
probably not always) a vertex operator superalgebra --- in particular
${}^1\widehat{{\cal V}}$ is an integral form of the Moonshine module 
$V^\natural$. Each ${}^g\widehat{{\cal V}}$ carries a representation of a
central extension of $C_{\M}(g)$ by $C_n$. Define the graded trace 
$$B(g,h;\tau)=q^{-1}\sum_{i\in {1\over n}\Z}{\rm ch}_{{}^g\widehat{{\cal V}}_i}
(h)\,q^i\ .$$
If $g,h\in\M$ commute and have
coprime orders, then $B(g,h;\tau)=T_{gh}(\tau)$. If all $q$-coefficients of
$T_g$ are nonnegative, then the `odd' part of  ${}^g\widehat{{\cal V}}$
vanishes, and ${}^g\widehat{{\cal V}}$ is the $g$-twisted
module $V^\natural(g)$ of} [{\bf 36}]. {\it If $g$ has prime order $p$, then the 
reduction mod $p$ of ${}^g\widehat{{\cal V}}$ is the modular vertex operator
superalgebra ${}^g{\cal V}$ of Theorem 7.}\medskip

When we say ${}^1\widehat{{\cal V}}$ is an integral form for $V^\natural$,
we mean that ${}^1\widehat{{\cal V}}$ has the same structure as a VOA, with
everything defined over $\Z$, and tensoring it with $\C$ recovers $V^\natural$.
This remarkable conjecture, which tries to explain Theorem 7, is completely
open.

\medskip{{\it 5.5.  The geometry of Moonshine.}}
Algebra is the mathematics of structure, and so of course it has a profound
relationship with every area of mathematics. Therefore the trick for finding
possible fingerprints of Moonshine in say geometry is to look there for modular
functions. And that search quickly leads to the elliptic genus.

For details see e.g.\ [{\bf 55}], [{\bf 108}], [{\bf 112}]. All manifolds here
are compact, oriented and differentiable. In Thom's cobordism
ring $\Omega$, elements are equivalence classes of cobordant manifolds, addition
is connected sum, and multiplication is Cartesian product.
The {\it universal elliptic genus} $\phi(M)$ is a ring homomorphism
from $\Q\otimes\Omega$ to the ring of power series in $q$, which sends $n$-dimensional
manifolds with spin connections to a weight $n/2$ modular form of $\Gamma_0(2)$ with 
integer coefficients. Several variations and generalisations have been
introduced, e.g.\ the Witten genus assigns spin manifolds with vanishing first
Pontryagin class a weight $n/2$ modular form of SL$_2(\Z)$ with integer
coefficients.

Several deep relationships between elliptic genera and the general material reviewed
elsewhere in this paper, have been uncovered. For instance, the important
rigidity property of the Witten genus with respect to any compact Lie group
action on the manifold,  is a consequence of the modularity of the characters
of affine algebras (our Theorem 4) [{\bf 81}]. The elliptic genus of a
manifold $M$ has been interpreted as the graded dimension of a vertex operator
superalgebra constructed from $M$ [{\bf 111}]. Seemingly related to this,
[{\bf 18}] recovered the elliptic genus of a Calabi--Yau manifold $X$ from the
sheaf of vertex algebras in the chiral de Rham complex [{\bf 85}] attached to
$X$. Unexpectedly, the elliptic genus of even-dimensional projective spaces
$P^{2n}$ has nonnegative coefficients and in fact equals the graded
dimension of a certain vertex algebra [{\bf 86}]; this suggests interesting
representation-theoretic questions in the spirit of Monstrous Moonshine.
In physics, elliptic
genera arise as partition functions of $N=2$ superconformal field theories
[{\bf 120}]. Mason's constructions [{\bf 88}] associated to Moonshine for the Mathieu
group $M_{24}$ have been interpreted as providing a geometric model (`elliptic
system') for elliptic
cohomology Ell${}^*(BM_{24})$ of the classifying space of $M_{24}$ [{\bf 112}],
[{\bf 39}].
The Witten genus (normalised by $\eta^8$) of the Milnor--Kervaire manifold
$M_0^8$, an 8-dimensional manifold built from the $E_8$ diagram, equals
$j^{{1\over 3}}$ [{\bf 55}] (recall (3.5)).

Hirzebruch's `prize question'
(p.86 of [{\bf 55}]) asks for the construction of a 24-dimensional manifold
$M$ with Witten genus $J$ (after being normalised by $\eta^{24}$). We would
like $\M$ to act on $M$ by diffeomorphisms, and the twisted Witten genera to be
the McKay--Thompson series $T_g$. It would also be nice to associate Norton's
series $N(g,h)$ to this Moonshine manifold. Constructing such a manifold
is perhaps the remaining Holy Grail of Monstrous Moonshine.

Hirzebruch's question was partially answered by Mahowald and Hopkins [{\bf 84}], 
who constructed a manifold with Witten genus $J$, but couldn't show it would
support an effective action of $\M$. Related work is [{\bf 3}], who constructed
several actions of $\M$ on e.g.\ 24-dimensional manifolds (but none of which
could have genus $J$), and [{\bf 73}], who showed the graded dimensions of the
subspaces $V^\natural_\pm$  of the Moonshine module are twisted $\widehat{A}$-genera
of Milnor--Kervaire's manifold $M^8_0$ (the $\widehat{A}$-genus is the specialisation
of elliptic genus to the cusp $\i\infty$).

\smallskip There has been a second conjectured relationship between geometry and
Monstrous Moonshine. {\it Mirror symmetry} says that most Calabi--Yau manifolds come
in closely related pairs. Consider a 1-parameter family $X_z$ of 
Calabi--Yau manifolds, with mirror $X^*$ given by the resolution of an orbifold
$X/G$ for $G$ finite and abelian. Then the Hodge numbers $h^{1,1}(X)$ and
$h^{2,1}(X^*)$ will be equal, and more precisely the moduli space of (complexified) K\"ahler
structures on $X$ will be locally isometric to the moduli space of complex
structures on $X^*$. The `mirror map' $z(q)$, which can be defined using
the Picard--Fuchs equation [{\bf 95}], gives a canonical map between those moduli
spaces. For example, $x_1^4+x_2^4+x_3^4+x_4^4+z^{-1/4}x_1x_2x_3x_4=0$
is such a family of K3 surfaces, where $G=C_4\times C_4$. Its mirror map
is given by
$$z(q)=q-104q^2+6444q^3-311744q^4+13018830q^5-493025760q^6+\cdots\ .\eqno(5.6)$$  

Lian--Yau [{\bf 80}] noticed that the reciprocal $1/z(q)$ of the mirror map
in (5.6) equals the McKay--Thompson series $T_g(\tau)+104$ for $g$ in class
2A of $\M$. After looking
at several other examples with similar conclusions, they proposed their
 {\it Mirror-Moonshine Conjecture}: The reciprocal $1/z$ of the mirror
map of a 1-parameter family of K3 surfaces with an orbifold mirror, will be
a McKay--Thompson series (up to an additive constant). 
 
A counterexample (and more examples) are given in \S7 of [{\bf 118}]. In 
particular, although there are relations between mirror symmetry and 
modular functions (see e.g.\ [{\bf 51}] and [{\bf 54}]), 
there doesn't seem to be any special relation with the Monster. Doran [{\bf
40}] `demystifies the Mirror-Moonshine phenomenon' by finding necessary and
sufficient conditions for $1/z$ to be a modular function for a modular
group commensurable with SL$_2(\Z)$. 

\bigskip\centerline{{\it 6. The physics of Moonshine}}\medskip

The physical side (perturbative string theory, or equivalently conformal
field theory) of Moonshine was noticed early on, and has profoundly influenced
the development of Moonshine and VOAs. This is a very rich subject, which
we can only superficially touch on. The book [{\bf 32}], with its extensive
bibliography, provides an introduction but will be difficult reading for many 
mathematicians (as will this section!) The treatment in [{\bf 45}] is more 
accessible and shows how naturally VOAs arise from the physics.
This effectiveness of physical interpretations isn't magic --- it merely
tells us that many of our finite-dimensional objects are seen much more
clearly when studied through infinite-dimensional structures (often by being
`looped'). Of course Moonshine, which teaches us to study the finite group
$\M$ via its infinite-dimensional module $V^\natural$, fits perfectly into
this picture.

A conformal field theory (CFT) is a quantum field theory on 2-dimensional
space-time, whose symmetries include the conformal transformations. In string
theory the basic objects are finite curves (`strings') rather than points
(`particles'), and the CFT lives on the surface traced by the strings as they
evolve (colliding and separating) through time. Each CFT is associated with
a pair ${\cal V}_L,{\cal V}_R$ of mutually commuting VOAs, called its {\it
chiral algebras} [{\bf 6}]. For example, strings living on a compact Lie group manifold
(the so-called Wess--Zumino--Witten model) will have chiral algebras given by
affine algebra VOAs. The space ${\cal H}$ of states for the
CFT carries a representation of ${\cal V}_L\otimes\overline{{\cal V}}_R$,
and many authors have (somewhat optimisticly) concluded that the study of
CFTs reduces to that of VOA representation theory. Rational VOAs correspond
to the important class of {\it rational} CFTs, where ${\cal H}$ decomposes into a
finite sum $\oplus M_L\otimes M_R$ of irreducible modules. The Virasoro
algebra (4.3) arises naturally in CFT through infinitesimal conformal transformations.
The vertex operator $Y(\phi,z)$, for the space-time parameter $z=e^{t+\i x}$,
is the quantum field which creates from the vacuum $|0\rangle\in{\cal H}$ the
state $|\phi\rangle\in{\cal H}$ at time $t=-\infty$: $|\phi\rangle={\rm lim}_{z
\rightarrow 0}Y(\phi,z)\,|0\rangle$. In particular, Borcherds' definition 
[{\bf 9}] of VOAs can be
interpreted as an axiomatisation of the notion of chiral algebra in CFT, and
for this reason alone is important.

In CFT, the Hauptmodul property of Moonshine is hard to interpret, and a less
direct formulation like that in [{\bf 116}] is needed. However, both the statement
and proof of Theorem 5 are natural from the CFT framework (see [{\bf 45}])
--- e.g.\ the modularity of the series $T_g$ and $N(g,h)$ are automatic in
CFT. This modularity arises in CFT through the equivalence of the Hamiltonian
formulation, which describes concretely the graded spaces we take traces on
(and hence the coefficients of our $q$-expansions), and the Feynman path
formalism, which interprets these graded traces as sections over moduli spaces
(and hence makes modularity manifest). Beautiful reviews are sketched in 
[{\bf 119}], [{\bf 120}].

Because $V^\natural$ is so mathematically special, it may be expected that it
corresponds to interesting physics. Certainly it has been the subject of
some speculation. There will be a $c=24$ rational CFT whose chiral algebra
${\cal V}_L$ and state space ${\cal H}$ are both  $V^\natural$, while
${\cal V}_R$ is trivial (this is possible
because $V^\natural$ is holomorphic). This CFT is nicely described in [{\bf 34}];
see also [{\bf 35}].
The Monster is the symmetry of that CFT, but the Bimonster $\M\wr C_2$ will
be the symmetry of a rational CFT with ${\cal H}=V^\natural\otimes \overline{V^\natural}$. 
The paper [{\bf 27}] finds a family of D-branes for the latter theory which
are in one-to-one correspondence with the elements of $\M$, and their `overlaps'
$\langle\!\langle g\|q^{{1\over 2}(L_0+\bar{L}_0-{c\over 24})}\|h\rangle\!\rangle$
equal the McKay--Thompson series $T_{g^{-1}h}$. However, we still lack any explanation
as to why a CFT involving $V^\natural$ should yield interesting {\it physics}.

Almost every facet of Moonshine finds a natural formulation in CFT, where it
often was discovered first. For example, the `No-Ghost' Theorem of
Brower--Goddard--Thorn
was used to great effect in [{\bf 11}] to understand the structure of the 
Monster Lie algebra ${\frak m}$. On a finite-dimensional manifold $M$, the 
index of the Dirac operator $D$ in the heat kernel interpretation is a path
integral in supersymmetric quantum mechanics, i.e\ an integral over the
free loop space ${\cal L}M=\{\gamma:S^1\rightarrow M\}$; the string theory
version of this is that the index of the Dirac operator on ${\cal L}M$ should
be an integral over ${\cal L}({\cal L}M)$, i.e.\ over smooth maps of tori into
$M$, and this is just the elliptic genus, and explains why it should be
modular. The orbifold construction of [{\bf 36}]
comes straight from CFT (although [{\bf 43}]'s construction of $V^\natural$
predates CFT orbifolds by a year and in fact influenced their development
in physics). That said, the translation process from
physics to mathematics of course is never easy --- Borcherds' definition
[{\bf 9}] is a prime example! 

But from this standpoint,
what is most exciting is what hasn't yet been fully exploited. String theory
tells us that CFT can live on any surface $\Sigma$. The VOAs, including the
geometric VOAs of [{\bf 59}], capture CFT in genus 0. The graded dimensions and
traces considered above concern CFT quantities (`conformal blocks') at genus
1: $\tau\mapsto e^{2\pi\i\tau}$ maps $\H$ onto a cylinder, and the trace
identifies the two ends. But there are analogues of all this at higher genus 
[{\bf 123}] (though the formulas can rapidly
become awkward). For example, the graded dimension of e.g.\ the $V^\natural$
CFT in genus 2 is computed in [{\bf 117}], and involves e.g.\ Siegel theta
functions. The orbifold theory in \S5.1 is genus 1: each `sector' $(g,h)$
corresponds to a homomorphism from the fundamental group $\Z^2$ of the torus
into the orbifold group $G$ (e.g\ $G=\M$) --- $g$ and $h$ are the targets of
the two generators of $\Z^2$ and hence must commute. More generally, the
sectors will correspond to each homomorphism $\varphi:\pi_1(\Sigma)\rightarrow
G$, and to each we will get a higher genus trace ${\cal Z}(\varphi)$, which
will be a function on the Teichm\"uller space $T_g$ (generalising the upper half-plane
${\Bbb H}$ for genus 1). The action of SL$_2(\Z)$ on the $N(g,h)$
generalises to the action of the mapping class group on $\pi_1(\Sigma)$ and
 $T_g$. See e.g.\ [{\bf 4}] for some thoughts in this direction.

\bigskip\centerline{{\it 7. Conclusion}}\medskip

There are different basic aspects to Monstrous Moonshine: (i) why modularity enters
at all; (ii) why in particular we have genus 0; and (iii) what does it have to 
do with the Monster. We understand (i) best. There will be a Moonshine-like 
relation between any (subgroup of the) automorphism
group of any rational VOA, and the characters $\chi_M$, and the
same can be expected to hold of the orbifold characters ${\cal Z}$ in \S5.1. 

To prove the genus 0 property of the $T_g$, we needed recursions obtained one way 
or another from the Monster Lie algebra ${\frak m}$, and from these we apply
Theorem 6. These recursions are very special, but so presumably is the
genus 0 property. The suggestion of [{\bf 20}] though is that we may be able
to considerably simplify this part of the argument.

Every group known to have rich Moonshine properties is contained in the Monster.
Our understanding of this seemingly central role of $\M$ is the poorest of those three
aspects.

It should be clear
from this review, of the central role VOAs play in our current understanding
of Moonshine. The excellent review [{\bf 39}] makes this point even more 
forcefully. It can be  (and has been) questioned
though whether the full and difficult machinery of VOAs is really needed to 
understand
this, i.e.\ whether we really have isolated the key conjunction of properties
needed for Moonshine to arise. CFT has been an invaluable guide thus far, but
perhaps we are a little too steeped in its lore.

Moonshine (in its more general sense) is a relation between algebra and number theory, and its impact on
algebra has been dramatic (e.g.\ VOAs, $V^\natural$, Borcherds--Kac--Moody algebras).
Its impact on number theory has been far less so. This may merely be a
temporary accident due to the backgrounds of most researchers (including the
mathematical physicists) working to date in the area. But  
the most exciting prospects for the future of Moonshine (in this writer's
opinion) are in the direction of 
number theory.
Hints of this future can be found in e.g.\ [{\bf 121}], [{\bf 41}], [{\bf 14}], 
[{\bf 33}], [{\bf 52}], [{\bf 94}].

\bigskip{{\it Acknowledgements.}} This paper was written at IHES
and at U of Wales Swansea, and I thank both for their hospitality. Much of
what I've learned of Moonshine came from conversations over the years with 
Chris Cummins and John McKay, and I warmly thank both. I appreciate the many
comments and suggestions provided by the readers of the first draft;
in this vein I would
particularly like to mention James Lepowsky and Fyodor Malikov.
My research is supported in part by NSERC.

\bigskip\centerline{{\it References}\footnote{$^\dagger$}{{\smal
For a more comprehensive Moonshine bibliography, especially for older
papers, see}}\footnote{}{{\smal 
http:$/\!/$cicma.mathstat.concordia.ca/faculty/cummins/moonshine.html}}}\medskip

\item{{\bf 1.}} {\smcap D.\ Alexander, C.\ Cummins, J.\ McKay,} and {\smcap
C.\ Simons}, `Completely replicable functions', {\it Groups, Combinatorics
and Geometry} (ed.\ M.W.\ Liebeck and J.\ Saxl, Cambridge Univ.\ Press, 1992)
87--98.

\item{{\bf 2.}} {\smcap O.\ Alvarez}, `Conformal anomalies and the index theorem',
{\it Nucl.\ Phys.} B286 (1987) 175--188.

\item{{\bf 3.}} {\smcap M.G.\ Aschbacher}, `Finite groups acting on homology 
manifolds',
{\it  Olga Taussky-Todd: in Memoriam, Pacific J.\ Math.\ Special Issue} (ed.\
M.\ Aschbacher et al, Pacific Journal of Mathematics, Berkeley, 1997) 3--36.

\item{{\bf 4.}} {\smcap P.\ B\'antay}, `Higher genus moonshine', {\it Moonshine,
the Monster, and Related Topics, Contemp.\
Math.} 193 (Amer.\ Math.\ Soc., Providence, 1996) 1--8.

\item{{\bf 5.}} {\smcap P.\ B\'antay}, `The kernel of the modular representation 
and the Galois action in RCFT', {\it Commun.\ Math.\ Phys.} 233 (2003) 423--438.

\item{{\bf 6.}} {\smcap A.\ Belavin, A.M.\ Polyakov,} and {\smcap A.B.\
Zamolodchikov}, `Infinite conformal symmetries in two-dimensional quantum
field theory', {\it Nucl.\ Phys.} {\bf B241} (1984) 333--380.

\item{{\bf 7.}} {\smcap J.S.\ Birman}, {\it Braids, Links and Mapping Class Groups}
(Princeton Univ.\ Press, 1974).

\item{{\bf 8.}} {\smcap S.\ Bloch}, `Zeta values and differential operators
on the circle', {\it J.\ Alg.} 182 (1996) 476--500.

\item{{\bf 9.}} {\smcap R.E.\ Borcherds}, `Vertex algebras, Kac--Moody algebras, and the
Monster', {\it Proc.\ Natl.\ Acad.\ Sci.\ (USA)} 83 (1986) 3068--3071.

\item{{\bf 10.}} {\smcap R.E.\ Borcherds}, `Generalized Kac--Moody Lie algebras',
{\it J.\ Alg.} 115 (1988) 501--512.

\item{{\bf 11.}} {\smcap R.E\ Borcherds}, `Monstrous moonshine and monstrous Lie superalgebras',
{\it Invent.\ Math.} 109 (1992) 405--444.

\item{{\bf 12.}} {\smcap R.E\ Borcherds}, `Sporadic groups and string theory',
{\it First European Congress of Math.\ Paris (1992)}, Vol.I (Birkh\"auser,
Basel, 1994) 411--421.

\item{{\bf 13.}} {\smcap R.E.\ Borcherds}, `Modular moonshine III', {\it Duke Math.\ J.}\ 93 (1998)
129--154.

\item{{\bf 14.}} {\smcap R.E.\ Borcherds}, `Automorphic forms with 
singularities on Grassmannians', {\it Invent.\ Math.} 132 (1998) 491--562.

\item{{\bf 15.}} {\smcap R.E.\ Borcherds}, `What is Moonshine?', {\it Proc. Intern.\ 
Congr.\ Math.\ 1998 (Berlin)}, Vol. 1 (Documenta Mathematica, Bielefeld, 1998)
607--615.

\item{{\bf 16.}} {\smcap R.E.\ Borcherds}, `Problems in Moonshine', 
{\it First Intern.\ Congr.\  Chinese Math.} 
(Amer. Math. Soc., Providence, 2001) 3--10. 
 
\item{{\bf 17.}} {\smcap R.E.\ Borcherds} and {\smcap A.J.E.\ Ryba},
`Modular moonshine II', {\it Duke Math.\ J.} 83 (1996) 435--459.

\item{{\bf 18.}} {\smcap L.A.\ Borisov} and {\smcap A.\ Libgober}, `Elliptic
genera of toric varieties and applications to mirror symmetry', {\it Invent.\
Math.} 140 (2000) 453--485.

\item{{\bf 19.}} {\smcap S.-P. Chan}, {\smcap M.-L. Lang}, and {\smcap C.-H.\ 
Lim}, `Some modular functions
associated to Lie algebra $E_8$', {\it Math.\ Z.} 211 (1992) 223--246.

\item{{\bf 20.}} {\smcap H.\ Cohn} and {\smcap J.\ McKay}, `Spontaneous
generation of modular invariants', {\it Math.\ of Comput.} 65 (1996) 1295--1309.

\item{{\bf 21.}} {\smcap J.H.\ Conway}, `Monsters and Moonshine', {\it Math.\ Intell.}\ 2 (1980) 165--171.

\item{{\bf 22.}} {\smcap J.H.\ Conway}, {\smcap R.T.\ Curtis}, {\smcap S.P.\ Norton},
{\smcap R.A.\ Parker} and {\smcap R.A.\
Wilson}, {\it An Atlas of Finite Groups} (Clarendon Press, Oxford, 1985).

\item{{\bf 23.}}  {\smcap J.H.\ Conway}, {\smcap J.\ McKay} and  {\smcap A.\ 
Sebbar}, `On the discrete groups of moonshine', {\it Proc.\ Amer.\ Math.\ Soc.}
(to appear).

\item{{\bf 24.}} {\smcap J.H.\ Conway} and {\smcap S.P.\ Norton}, `Monstrous Moonshine', {\it Bull.\ London
Math.\ Soc.} 11 (1979) 308--339.

\item{{\bf 25.}} {\smcap J.H.\ Conway}, {\smcap S.P.\ Norton}, and {\smcap
L.H.\  Soicher}, `The Bimonster, the group $Y_{555}$, and the projective plane 
of order $3$', {\it Computers in Algebra} (Dekker, New York, 1988) 27--50.

\item{{\bf 26.}} {\smcap J.H.\ Conway} and {\smcap N.J.A.\ Sloane}, {\it Sphere Packings,
Lattices and Groups}, 3rd edn (Springer, Berlin, 1999).

\item{{\bf 27.}} {\smcap B.\ Craps}, {\smcap M.R.\ Gaberdiel}, {\smcap
J.A.\ Harvey}, `Monstrous branes', {\it Commun.\ Math.\ Phys.} 234 (2003)
229--251.

\item{{\bf 28.}} {\smcap C.J.\ Cummins}, `Congruence subgroups of groups
commensurable with PSL$(2,\Z)$ of genus 0 and 1', Preprint.

\item{{\bf 29.}} {\smcap C.J. Cummins} and {\smcap T.\ Gannon}, `Modular equations and the genus
zero property of moonshine functions', {\it Invent.\ Math.} 129 (1997) 413--443.

\item{{\bf 30.}} {\smcap C.J.\ Cummins} and {\smcap S.P.\ Norton}, `Rational
Hauptmodul are replicable', {\it Canad.\ J.\ Math.} 47 (1995) 1201--1218.

\item{{\bf 31.}} {\smcap C.W.\ Curtis} and {\smcap I.\ Reiner}, {\it Methods of
Representation Theory with Applications to Finite Groups and Orders}, Vol.I
(Wiley, New York, 1981). 

\item{{\bf 32.}} {\smcap P.\ Di Francesco}, {\smcap P.\ Mathieu} and {\smcap
D.\ S\'en\'echal}, {\it Conformal Field Theory} (Springer, New York, 1997).

\item{{\bf 33.}} {\smcap R.\ Dijkgraaf}, `The mathematics of fivebranes',
{\it Proc.\ Intern.\ Congr.\ Math.\ 1998 (Berlin)}, Vol.III (Documenta Mathematica,
Bielefeld, 1998) 133--142.

\item{{\bf 34.}} {\smcap L.\ Dixon}, {\smcap P.\ Ginsparg} and {\smcap J.A.\ Harvey},
`Beauty and the beast:
superconformal symmetry in a monster module', {\it Commun.\ Math.\ Phys.} 119 (1988)
221--241.

\item{{\bf 35.}} {\smcap L.\ Dolan}, {\smcap P.\ Goddard} and {\smcap P.\
Montague}, `Conformal field theory of twisted vertex operators', {\it Nucl.\
Phys.} B338 (1990) 529--601.

\item{{\bf 36.}} {\smcap C.\ Dong}, {\smcap H.\ Li} and {\smcap G.\ Mason}, 
`Modular invariance of trace functions in orbifold theory and generalized 
moonshine', {\it Commun.\ Math.\ Phys.} 214 (2000) 1--56.


\item{{\bf 37.}} {\smcap C.\ Dong} and {\smcap G.\ Mason}, `Nonabelian
orbifolds and the boson-fermion correspondence', {\it Commun.\ Math.\ Phys.}
163 (1994) 523--559.

\item{{\bf 38.}} {\smcap C.\ Dong} and {\smcap G.\ Mason}, `An orbifold
theory of genus zero associated to the sporadic group $M_{24}$', {\it
Commun.\ Math.\ Phys.} 164 (1994) 87--104.

\item{{\bf 39.}} {\smcap C.\ Dong} and {\smcap G.\ Mason}, `Vertex operator algebras 
and moonshine: A survey', 
{\it Progress in Algebraic Combinatorics, Adv.\ Stud.\ Pure Math.}\ 24 (Math.\
Soc.\ Japan, Tokyo, 1996) 101--136.

\item{{\bf 40.}} {\smcap C.F.\ Doran}, `Picard--Fuchs uniformization
and modularity of the mirror map', {\it Commun.\ Math.\ Phys.} 212 (2000) 625--647.

\item{{\bf 41.}} {\smcap V.G.\ Drinfeld}, `On quasitriangular quasi-Hopf algebras
on a group that is closely related with Gal$\,\overline{\Q}/\Q$', {\it Leningrad.\
Math.\ J.} 2 (1991) 829--860.

\item{{\bf 42.}} {\smcap I.\ Frenkel}, {\smcap Y.-Z.\ Huang} and {\smcap J.\
Lepowsky}, `On Axiomatic Approaches
to Vertex Operator Algebras and Modules', {\it Mem.\ Amer.\ Math.\ Soc.} 104
(1993) 1--63.

\item{{\bf 43.}} {\smcap I.\ Frenkel}, {\smcap J.\ Lepowsky} and {\smcap A.\ Meurman}, `A natural representation
of the Fischer--Griess monster with the modular function $J$ as character',
{\it Proc.\ Natl.\ Acad.\ Sci.\ USA} 81 (1984) 3256--3260.

\item{{\bf 44.}} {\smcap I.\ Frenkel}, {\smcap J.\ Lepowsky} and {\smcap A.\ Meurman},
{\it Vertex Operator Algebras and the Monster} (Academic Press, New York, 1988).

\item{{\bf 45.}} {\smcap M.R.\ Gaberdiel} and {\smcap P.\ Goddard}, `Axiomatic
conformal field theory', {\it Commun.\ Math.\ Phys.} 209 (2000) 549--594.

\item{{\bf 46.}} {\smcap T.\ Gannon}, {\it Moonshine Beyond the Monster}
(Cambridge Univ.\ Press, to appear).

\item{{\bf 47.}} {\smcap G.\ Glauberman} and {\smcap S.P.\ Norton}, `On McKay's connection
between the affine $E_8$ diagram and the Monster', {\it Proc.\ on Moonshine and
Related Topics} (Amer.\ Math.\ Soc., Providence, 2001) 37--42.

\item{{\bf 48.}} {\smcap P.\ Goddard}, `The work of Richard Ewen Borcherds', 
{\it Proc.\ Intern.\ Congr.\ Math.\ 1998 (Berlin)}, Vol.\ I 
(Documenta Mathematica, Bielefeld, 1998) 99--108.

\item{{\bf 49.}} {\smcap D.\ Gorenstein}, {\it Finite Simple Groups: An Introduction to their
Classification} (Plenum, New York, 1982).

\item{{\bf 50.}} {\smcap R.\ Griess}, `The friendly giant', {\it Invent.\ Math.} 
68 (1982) 1--102.

\item{{\bf 51.}} {\smcap V.A.\ Gritsenko} and {\smcap V.V.\ Nikulin}, `The
arithmetic mirror symmetry and Calabi--Yau manifolds', {\it Commun.\ Math.\
Phys.} 210 (2000) 1--11. 

\item{{\bf 52.}} {\smcap S.\ Gukov} and {\smcap C.\ Vafa}, `Rational conformal
field theories and complex multiplication', Preprint (arXiv: hep-th/0203213).

\item{{\bf 53.}} {\smcap K.\ Harada}, {\smcap M.\ Miyamoto}, and {\smcap H.\
Yamada}, `A generalization of Kac--Moody algebras', {\it Groups, Difference Sets,
and the Monster} (de Gruyter, Berlin, 1996) 377--408.

\item{{\bf 54.}} {\smcap J.A.\ Harvey} and {\smcap G.\ Moore}, `Algebras, BPS
states, and strings', {\it Nucl.\ Phys.} B463 (1996) 315--368.

\item{{\bf 55.}} {\smcap F.\ Hirzebruch}, {\smcap T.\ Berger}, and {\smcap R.\ Jung}, {\it Manifolds and Modular
Forms} (Friedr.\ Vieweg \& Sohn, 1991).


\item{{\bf 56.}} {\smcap G.\ H\"ohn}, `The group of symmetries of the shorter
Moonshine module', Preprint (arXiv:math.QA/0210076).

\item{{\bf 57.}} {\smcap T.\ Hsu}, {\it Quilts: Central Extensions, Braid Actions, and Finite
Groups}, {\it Lecture Notes in Mathematics 1731} (Springer, Berlin, 2000).

\item{{\bf 58.}} {\smcap Y.-Z.\ Huang}, `Binary trees and finite-dimensional
Lie algebras', {\it Algebraic Groups and their Generalizations: Quantum and
Infinite-Dimensional Methods}, Proc.\ Symp.\ Pure Math.\ 56, pt.\ 2 (Amer.\
Math.\ Soc., Providence, 1994) 337--348. 

\item{{\bf 59.}} {\smcap Y.-Z.\ Huang}, {\it Two-dimensional Conformal Geometry
 and Vertex Operator Algebras} (Birkh\"auser, Boston, 1997).

\item{{\bf 60.}} {\smcap A.A.\ Ivanov}, `Geometric presentations of groups with an application
to the Monster', {\it Proc.\ Intern.\ Congr.\ Math.\ 1990 (Kyoto)}, Vol.\
II (Springer, Hong Kong, 1991) 1443--1453.

\item{{\bf 61.}} {\smcap A.A.\ Ivanov},  `$Y$-groups via transitive extension', 
{\it J.\ Alg.} 218 (1999) 412--435.

\item{{\bf 62.}} {\smcap R.\ Ivanov} and {\smcap M.P.\ Tuite}, `Some irrational
generalized Moonshine from orbifolds', {\it Nucl.\ Phys.} B635 (2002) 473--491.

\item{{\bf 63.}} {\smcap E.\ Jurisich}, `An exposition of generalized Kac--Moody
algebras', {\it Contemp.\ Math.} 194 (Amer.\ Math.\ Soc., Providence, 1996)
121--159.

\item{{\bf 64.}} {\smcap E.\ Jurisich}, {\smcap J.\ Lepowsky} and {\smcap 
R.L.\ Wilson}, `Realizations of the Monster Lie algebra', {\it Selecta Math.\
(NS)} 1 (1995) 129--161.

\item{{\bf 65.}} {\smcap V.G.\ Kac}, `An elucidation of: Infinite-dimensional algebras,
Dedekind's $\eta$-function, classical M\"obius function and the very strange
formula. $E_8^{(1)}$ and the cube root of the modular invariant $j$',
{\it Adv.\ in Math.} 35 (1980) 264--273.

\item{{\bf 66.}} {\smcap V.G.\ Kac}, {\it Infinite Dimensional Lie Algebras}, 3rd edn, 
(Cambridge Univ.\ Press, 1990).

\item{{\bf 67.}} {\smcap V.G.\ Kac}, {\it Vertex Operators for Beginners}, Univ.\
Lecture Series, Vol.\ 10 (Amer.\ Math.\ Soc., Providence, 1997).

\item{{\bf 68.}} {\smcap V.G.\ Kac} and {\smcap S.-J.\ Kang}, `Trace formula for
graded Lie algebras and Monstrous Moonshine', {\it Representations of Groups, Banff
(1994)} (Amer.\ Math.\ Soc., Providence, 1995) 141--154.

\item{{\bf 69.}} {\smcap V.G.\ Kac} and {\smcap D.H.\ Peterson}, `Infinite dimensional 
Lie algebras, theta functions, and modular forms', {\it Adv.\ in
Math.} 53 (1984) 125--264.

\item{{\bf 70.}} {\smcap S.\ Kass}, {\smcap R.V.\ Moody}, {\smcap J.\ Patera} and
{\smcap  R.\ Slansky}, {\it Affine Lie Algebras,
Weight Multiplicities, and Branching Rules}, Vol.\ 1 (Univ.\  California
Press, Berkeley, 1990).

\item{{\bf 71.}} {\smcap M.\ Knopp} and {\smcap G.\ Mason}, `Generalized
modular forms', {\it J.\ Number Theory} 99 (2003) 1--28.

\item{{\bf 72.}} {\smcap D.N.\ Kozlov}, `On completely replicable functions and extremal
poset theory', MSc thesis, Univ.\ of Lund, Sweden, 1994.

\item{{\bf 73.}} {\smcap R.\ Kultze}, `Elliptic genera and the moonshine module',
{\it Math.\ Z.} 223 (1996) 463--471.

\item{{\bf 74.}} {\smcap M.-L.\ Lang}, `On a question raised by Conway--Norton', 
{\it J. Math.\ Soc.\ Japan} 41 (1989) 263--284.

\item{{\bf 75.}} {\smcap S.\ Lang}, {\it Elliptic Functions}, 2nd edn 
(Springer, New York, 1997).

\item{{\bf 76.}} {\smcap J.\ Lepowsky}, `Euclidean Lie algebras and the modular 
function $j$', {\it The Santa Cruz Conference on Finite Groups},
{\it Proc.\ Sympos.\ Pure Math.} 37 (Amer.\ Math.\ Soc., Providence, 1980)
567--570.

\item{{\bf 77.}} {\smcap J.\ Lepowsky}, `Vertex operator algebras and the zeta
function', {\it Recent Developments in Quantum Affine Algebras and Related Topics, 
Contemp.\ Math.} 248 (Amer.\ Math.\ Soc., Providence, 1999) 327--340.

\item{{\bf 78.}} {\smcap J.\ Lepowsky}, `The work of Richard E.\ Borcherds',
 {\it Notices Amer.\ Math.\ Soc.} 46 (1999) 17--19.

\item{{\bf 79.}} {\smcap J.\ Lepowsky} and {\smcap H.\ Li}, {\it Introduction
to Vertex Operator Algebras and their Representations} (Birkh\"auser, Boston,
2004).

\item{{\bf 80.}} {\smcap B.H.\ Lian} and {\smcap S.-T.\ Yau}, `Arithmetic
properties of mirror map and quantum coupling', {\it Commun.\ Math.\ Phys.}
176 (1996) 163--191.

\item{{\bf 81.}} {\smcap B.H.\ Lian} and {\smcap G.J.\ Zuckerman}, `Moonshine
cohomology', 
{\it Moonshine and Vertex Operator Algebras}  
(Surikaisekikenkyusho Kokyuroku  No.\ 904, Kyoto, 1995) 87--115. 
 
\item{{\bf 82.}} {\smcap K.\ Liu}, `On modular invariance and rigidity theorems',
{\it J.\ Diff.\ Geom.} 41 (1995) 343--396.

\item{{\bf 83.}} {\smcap I.G.\ Macdonald}, `Affine root systems and Dedekind's
$\eta$-function', {\it Invent.\ Math.} 15 (1972) 91--143.

\item{{\bf 84.}} {\smcap M.\ Mahowald} and {\smcap M.\ Hopkins},
`The structure of 24 dimensional manifolds having normal bundles which lift to $B{\rm O}[8]$',
{\it Recent Progress in Homotopy Theory, Contemp.\ Math.} 293 
(Amer.\ Math.\ Soc., Providence, 2002) 89--110.

\item{{\bf 85.}} {\smcap F.\ Malikov, V.\ Schechtman} and {\smcap A.\ Vaintrob},
`Chiral de Rham complex', {\it Commun.\ Math.\ Phys.} {\bf 204} (1999) 439--473.

\item{{\bf 86.}} {\smcap F.\ Malikov} and {\smcap V.\ Schechtman}, `Deformations
of vertex algebras, quantum cohomology of toric varieties, and elliptic
genus', {\it Commun.\ Math.\ Phys.} 234 (2003) 77--100.

\item{{\bf 87.}} {\smcap G.\ Mason}, `Finite groups and modular functions', {\it 
The Arcata Conference on Representations of Finite Groups, Proc.\ Sympos.\
 Pure Math.} 47 (Amer.\ Math.\ Soc., Providence, 1987) 181--209.

\item{{\bf 88.}} {\smcap G.\ Mason}, 
`On a system of elliptic modular forms attached to the large Mathieu group',
{\it Nagoya Math. J.} 118 (1990) 177--193.

\item{{\bf 89.}} {\smcap J.\ McKay}, `Graphs, singularities, and finite groups',
{\it The Santa Cruz Conference on Finite Groups,
Proc.\ Sympos.\ Pure Math.} 37 (Amer.\ Math.\ Soc., Providence, 1980) 183--186.

\item{{\bf 90.}} {\smcap J.\ McKay}, `The essentials of Monstrous Moonshine', {\it Groups and
Combinatorics -- in Memory of M.\ Suzuki, Adv.\ Studies
in Pure Math.} 32 (Math.\ Soc.\ Japan, Tokyo, 2001) 347--353.

\item{{\bf 91.}}  {\smcap J.\ McKay} and {\smcap H.\ Strauss}, `The q-series
of monstrous moonshine and the decomposition of the head characters', {\it
Commun.\ Alg.} 18 (1990) 253--278.

\item{{\bf 92.}} {\smcap H.\ McKean} and {\smcap V.\ Moll}, {\it Elliptic Curves:
Function Theory, Geometry, Arithmetic} (Cambridge Univ.\ Press, 
1999).

\item{{\bf 93.}} {\smcap M.\ Miyamoto}, `21 involutions acting on the Moonshine
module', {\it J.\ Alg.} 175 (1995) 941--965.

\item{{\bf 94.}} {\smcap G.W.\ Moore}, `Les Houches lectures on strings and
arithmetic', Preprint (arXiv: hep-th/0401049).

\item{{\bf 95.}} {\smcap D.R.\ Morrison}, `Picard--Fuchs equations and mirror
maps for hypersurfaces', {\it Essays on Mirror Manifolds} (ed.\ S.-T.\ Yau, 
Intern.\ Press, Hong Kong, 1992) 241--264.

\item{{\bf 96.}} {\smcap S.P.\ Norton}, `More on Moonshine', In: {\it Computational
Group Theory} (ed.\ M.D.\ Atkinson, Academic Press, New York, 1984) 185--193.

\item{{\bf 97.}} {\smcap S.P.\ Norton}, `Generalized moonshine', {\it The Arcata Conference on
Representations of Finite Groups, Proc.\ Sympos.\ Pure Math.}
47 (Amer.\ Math.\ Soc., Providence, 1987) 208--209.

\item{{\bf 98.}} {\smcap S.P.\ Norton}, `Constructing the Monster', {\it
Groups, Combinatorics, and Geometry} (ed.\ M.W.\ Liebeck and J.\ Saxl,
Cambridge Univ.\ Press, 1992) 63--76.

\item{{\bf 99.}} {\smcap S.P.\ Norton}, `From moonshine to the Monster', {\it Proc.\
on Moonshine and Related Topics} (Amer.\ Math.\ Soc, Providence, 2001)
 163--171.

\item{{\bf 100.}} {\smcap A.\ Ogg}, {\it Modular forms and Dirichlet series}
 (W. A. Benjamin, New York, 1969).

\item{{\bf 101.}} {\smcap L.\ Queen}, `Modular functions arising from some finite groups',
{\it Math.\ of Comput.} 37 (1981) 547--580.

\item{{\bf 102.}} {\smcap U.\ Ray}, `Generalized Kac--Moody algebras
and some related topics', {\it Bull.\ Amer.\ Math.\ Soc.} 38 (2000) 1--42.

\item{{\bf 103.}} {\smcap W.F.\ Reynolds}, `Thompson's characterization of 
characters and sets of primes', {\it J.\ Alg.} 156 (1993) 237--243.

\item{{\bf 104.}} {\smcap S.-S.\ Roan}, `Mirror symmetry of elliptic curves
and Ising model', {\it J.\ Geom.\ Phys.} 20 (1996) 273--296.

\item{{\bf 105.}} {\smcap A.J.E.\ Ryba}, `Modular moonshine?', {\it Moonshine, the Monster,
and Related Topics, Contemp.\ Math.} 193 (Amer.\ Math.\ Soc., Providence, 1996)
307--336.

\item{{\bf 106.}} {\smcap A.N.\ Schellekens}, `Meromorphic $c=24$ conformal
field theories', {\it Commun.\ Math.\ Phys.} 153 (1993) 159--185.

\item{{\bf 107.}} {\smcap G.\ Segal}, `The definition of conformal field theory',
{\it Differential Geometric Methods in Theoretical Physics} (ed.\ K.\ Bleuler and
M.\ Werner, Academic Press, Boston, 1988) 165--171.

\item{{\bf 108.}} {\smcap G.\ Segal}, `Elliptic cohomology', {\it S\'eminaire
Bourbaki 1987-88, no.\ 695, Ast\'erisque} 161-162 (1988) 187--201.

\item{{\bf 109.}} {\smcap G.W.\ Smith}, `Replicant powers for higher genera',
{\it Moonshine, the Monster, and Related Topics, Contemp.\ Math.}
193 (Amer.\ Math.\ Soc., Providence, 1996) 337--352.

\item{{\bf 110.}} {\smcap S.D.\ Smith}, `On the head characters of the monster 
simple group',
{\it Finite Groups -- Coming of Age} (ed.\ J.\ McKay, Amer.\ Math.\ Soc., 
Providence, 1985) 303--313.

\item{{\bf 111.}} {\smcap H.\ Tamanoi}, {\it Elliptic Genera and Vertex Operator
Superalgebras, Lecture Notes in Mathematics 1704} (Springer, Berlin, 1999).

\item{{\bf 112.}} {\smcap C.B.\ Thomas}, {\it Elliptic Cohomology} (Kluwer, New York, 1999).

\item{{\bf 113.}} {\smcap J.G.\ Thompson}, `Finite groups and modular functions',
{\it Bull.\ London Math.\ Soc.} 11 (1979) 347--351.

\item{{\bf 114.}} {\smcap J.G.\ Thompson}, `Some numerology between the Fischer--Griess
Monster and the elliptic modular function', {\it Bull.\ Lond.\ Math.\ Soc.}
11 (1979) 352--353.

\item{{\bf 115.}} {\smcap J.G.\ Thompson}, `A finiteness theorem for subgroups
of PSL(2,$\R$) which are commensurable with PSL($2,\Z)$', {\it Santa Cruz 
Conference on Finite Groups,  Proc.\ Symp.\ Pure Math.} 37 (Amer.\ Math.\ 
Soc., Providence, 1980) 533--555.

\item{{\bf 116.}} {\smcap M.P.\ Tuite}, `On the relationship between Monstrous Moonshine
and the uniqueness of the Moonshine module', {\it Commun.\ Math.\ Phys.} 166 (1995)
495--532.

\item{{\bf 117.}} {\smcap M.P.\ Tuite}, `Genus two meromorphic conformal field
theories', {\it Proc.\ on Moonshine and Related Topics} (Amer.\ Math.\ Soc.,
Providence, 2001) 231--251.

\item{{\bf 118.}} {\smcap H.\ Verrill} and {\smcap N.\ Yui}, `Thompson series
and the mirror maps of pencils of K3 surfaces', {\it The Arithmetic and Geometry
of Algebraic Cycles} (ed.\ B.\ Gordon et al, Centre Res.\ Math.\ Proc.\ and
Lecture notes 24, 2000) 399--432.

\item{{\bf 119.}} {\smcap E.\ Witten}, `Physics and geometry', {\it Proc.\
Intern.\ Congr.\ Math.\ 1986 (Berkeley)} (Amer.\ Math.\ Soc., Providence, 1987).

\item{{\bf 120.}} {\smcap E.\ Witten}, `Elliptic genera and quantum field theory',
{\it Commun.\ Math.\ Phys.} 109 (1987) 525--536.

\item{{\bf 121.}} {\smcap E.\ Witten}, `Quantum field theory, Grassmannians,
and algebraic curves', {\it Commun.\ Math.\ Phys.} 113 (1988) 529--600.

 \item{{\bf 122.}} {\smcap E.\ Witten}, `Geometry and quantum field theory',
{\it Mathematics Into the Twenty-first Century}, Vol.II (Amer.\ Math.\ Soc.,
Providence, 1992) 479--491.

\item{{\bf 123.}} {\smcap Y.\ Zhu}, `Global vertex operators on Riemann 
surfaces', {\it Commun.\ Math.\ Phys.} 165 (1994) 485--531.

\item{{\bf 124.}} {\smcap Y.\ Zhu}, `Modular invariance of characters of vertex operator
algebras', {\it J.\ Amer.\ Math.\ Soc.} 9 (1996) 237--302.

\bigskip\bigskip{\it Terry Gannon}

{\it Department of Mathematics}

{\it University of Alberta}

{\it Edmonton, CANADA T6G 2G1}

\medskip tgannon@math.ualberta.ca

\end